\documentclass[10pt,reqno]{amsart}

\usepackage[text={160mm,240mm},centering]{geometry}            
\usepackage{comment}
\usepackage{graphicx}
\usepackage{amssymb}
\usepackage{epstopdf}
\usepackage{dsfont}
\usepackage{tikz-cd}
\usepackage{amsmath,amsthm}
\usepackage{amsmath}
\usepackage{mathabx}
\usepackage[mathscr]{eucal}
\usepackage[all]{xy}
\usepackage{mathrsfs}
\usepackage{hyperref}
\usepackage{color}
\definecolor{ccr}{RGB}{10,110,150}
\usepackage{hyperref}
\hypersetup{hypertex=true,colorlinks=true,linkcolor=ccr,anchorcolor=ccr,citecolor=ccr}

\usepackage{ulem}

\newcommand{\m}{\mathrm}
\newtheorem{theorem}{Theorem}[section]
\newtheorem{lemma}[theorem]{Lemma}
\numberwithin{equation}{section}

\newtheorem{corollary}[theorem]{Corollary}
\newtheorem{proposition}[theorem]{Proposition}

\newtheorem{question}[theorem]{Question}
\newtheorem{definition-lemma}[theorem]{Definition}
\newtheorem{definition-theorem}[theorem]{Definition}

\theoremstyle{definition}
\newtheorem{example}[theorem]{Example}
\newtheorem{definition}[theorem]{Definition}

\usepackage{bm}

\newtheorem{remark}[theorem]{Remark}

\newcommand{\fl}{\mathcal{FL}}

\newcommand{\xkn}{X_n(K),K_n}
\newcommand{\spc}{\mathrm{Spin}^c}
\newcommand{\rspc}{\underline{\spc}}
\newcommand{\hfk}{\mathrm{HFK}}
\newcommand{\hfl}{\mathrm{HFL}}
\newcommand{\cfk}{\mathrm{CFK}}
\allowdisplaybreaks

\title{On the nonvanishing of Floer lasagna modules}

\author{Kai Zhu}
\address{Beijing International Center for Mathematical Research, Peking University, Beijing 100871,  P. R. China.}
\email{kzhu25@stu.pku.edu.cn}

\begin{document}

	\maketitle
	
	\renewcommand{\abstractname}{Abstract}
	\begin{abstract}
		In this paper, we establish several nonvanishing results for Floer lasagna modules $\mathcal{FL}(W,K)$ for a smooth compact connected oriented $4$-manifold $W$ and a knot $K\subset \partial W$, via the evaluation maps to link Floer homology. We first establish a general criterion: if a certain Heegaard Floer cobordism map associated with a smooth $4$-manifold $W$ is nonzero, then $\mathcal{FL}(W,U)\neq 0$, where $U\subset \partial W$ is a local unknot. Applications include the nonvanishing of Floer lasagna modules for $\#^m\overline{\mathbb{C}P^2}$ and for knot traces $X_n(K)$ when $K$ is slice and $n$ is a non-positive integer. In another direction, we show that for every nontrivial knot \(K\) and every integer \(n\), the Floer lasagna module $\mathcal{FL}(X_n(K),K_n)$ is nonzero, where \(K_n\subset S_n^3(K)\) is the dual knot. Moreover, if \(K\) is an L-space knot, then $\mathcal{FL}(X_n(K),K_n)$ is infinite-dimensional for all integers $n\ge \max\{2g(K),1\}$.
	\end{abstract}
	
	\parskip=.3em
	
	\setcounter{page}{1}
	\section{Introduction}\label{sec1}

	In \cite{MWW22}, Morrison, Walker and Wedrich introduced the notion of skein
	lasagna module $\mathcal{S}^2_0(W,L)$, where $W$ is a smooth $4$-manifold and
	$L\subset \partial W$ is a boundary link.\footnote{Throughout the paper, by a
		smooth $4$-manifold we mean a compact, connected, oriented smooth
		$4$-manifold.} It is an invariant of the pair $(W,L)$ constructed from
	Khovanov--Rozansky homology. Later in \cite{RW26}, Ren and Willis defined a lasagna
	$s$-invariant from $\mathcal{S}^2_0(W,L)$ and gave the first analysis-free
	proof of the existence of exotic $4$-manifolds. In principle, one can use
	different TQFTs to define various lasagna modules. A natural analogue is to use
	the link Floer TQFT defined by Zemke in \cite{Zem19a,Zem19b}, the corresponding
	Floer lasagna module $\fl(W,L)$ was first introduced by Chen in
	\cite{chen2022floer}. When the boundary link is empty, we write
	$\fl(W)=\fl(W,\varnothing)$.
	
Informally, a lasagna filling of $(W,L)$ consists of a finite collection of
properly embedded surfaces in $W$ with the interiors of finitely many disjoint
input $4$-balls removed, such that the outer boundary is $L$ and the inner
boundary links are labelled by elements of the chosen link homology theory. The lasagna
module is generated by such fillings, modulo multilinearity in the labels and
local replacements.

	There are many computational results for the skein lasagna module, for example \cite{MN22,MWW23,RW26,SZ24}, since Khovanov-Rozansky homology is purely algebraic and combinatorial. However, the link Floer analogue is harder to compute, since it requires computing a lot of holomorphic data and is therefore more technically involved. At the time of writing, explicit computations of Floer lasagna modules remain limited. The following results are known:
	\begin{itemize}
		\item In \cite{chen2022floer}, Chen showed that
		$\fl(S^2\times \mathbb{D}^2,L)=0$ whenever the
		surgery-diagram representative of
		$L\subset S^1\times S^2$ intersects a spanning disk
		of the $0$-framed surgery unknot geometrically once.
		
		\item In \cite{McC26}, McCulloch showed that
		$\fl(X_{-n}(U))$ is infinite dimensional for every
		integer $n>0$, where $U$ denotes the unknot and
		$X_{-n}(U)$ denotes its $(-n)$-framed knot trace.
		
		\item More recently, Alishahi, Gorsky, and Liu
		\cite{AGL26}
		computed the Floer lasagna modules of
		$S^2\times\mathbb{D}^2$ and $X_{-1}(U)$ with empty
		boundary link in all homology classes and Alexander
		and Maslov gradings. More generally, they obtained
		a complete graded description of
		$\fl(X_{-n}(m(K)))$ for any plumbed
		L-space knot $K$ and integer
		$n>\max\{0,2g(K)-1\}$, where $m(K)$ denotes
		the mirror of $K$.
	\end{itemize}

In this paper, we mainly use the evaluation map
\[
\m{ev}_{\mathfrak{s}}:\mathcal{FL}(W,L)\to \widehat{\hfl}(\partial W,L)
\]
for each $\spc$ structure $\mathfrak{s}$ on $W$, which generalizes the
evaluation map defined in \cite{chen2022floer}. Together with some computational
results for link Floer cobordism maps, this map allows us to obtain lower bounds
on the dimension of the Floer lasagna module and hence to prove several
nonvanishing results for Floer lasagna modules. More precisely, we construct
specific Floer lasagna fillings and show that their images under the evaluation
map are nonzero.

	For example, let \(W\) be a smooth  \(4\)-manifold with nonempty boundary \(Y\). Since $\partial W\ne\varnothing$, we may assume that \(W\) admits a handle decomposition with a single \(0\)-handle and no \(4\)-handle. Take the input ball \(B\) to be the single \(0\)-handle, and take a local unknot \(U\) away from the attaching spheres, labelled by \(1\in \widehat{\hfl}(S^3,U)\). The track of this local unknot \(U\) forms a product annulus \(U\times I\) from \((S^3,U)\) to \((Y,U)\). Hence we obtain a canonical Floer lasagna filling of \((W,U)\), which we denote by \(\mathcal{F}_{W,U}\). If \(Y=S^3\), we can further cap \(W\) by a \(4\)-handle and cap \(U\times I\) by a disk. This yields a canonical Floer lasagna filling of a closed manifold with empty boundary link. Denote the capped-off manifold by \(X\) and the resulting canonical Floer lasagna filling by \(\mathcal{F}_X\).
	
We use the evaluation map to show that the classes represented by the canonical Floer lasagna fillings $\mathcal{F}_{W,U}$ and $\mathcal{F}_X$, denoted by $[\mathcal{F}_{W,U}]$ and $[\mathcal{F}_X]$ respectively, are nonzero under certain conditions. In particular, we have the following criterion.
	
	\begin{theorem}\label{thmgeneralW}
		Let \(W\) be a smooth  \(4\)-manifold with connected boundary \(Y=\partial W\). Removing the interior of a \(4\)-ball \(B\subset W^\circ\) gives a cobordism \(W\setminus B^\circ:S^3\to Y\). Assume that there exists a \(\m{Spin}^c\) structure $\mathfrak{s}$ on $W$ such that the Heegaard Floer cobordism map
		\begin{equation}
			\widehat{F}_{W\setminus B^\circ,\mathfrak{s}}:\widehat{\m{HF}}(S^3)\to \widehat{\m{HF}}(Y)
		\end{equation}
		is nonzero. Then \(\left[\mathcal{F}_{W,U}\right]\ne 0\), where \(U\subset Y=\partial W\) is a local unknot. Hence \(\fl(W,U)\ne 0\).
		
		If \(Y=S^3\) and \(X=W\cup_{Y}B^4\) is the capped-off closed \(4\)-manifold obtained from \(W\), then \(\left[\mathcal{F}_X\right]\ne 0\), hence \(\fl(X)\ne 0\).
	\end{theorem}

	Combining various nonvanishing results for Heegaard Floer cobordism maps, we obtain the following,
	
	\begin{proposition}\label{thmcp2} 
		\label{thmstein}
		For any Stein domain \(X\), local unknot \(U\subset \partial X\), and integer \(n\ge 0\), we have
		\[
		\left[\mathcal{F}_{X\#^n\overline{\mathbb{C}P^2},U}\right]\ne 0,
		\quad\text{and}\quad
		\left[\mathcal{F}_{\#^n\overline{\mathbb{C}P^2}}\right]\ne 0.
		\]Hence $\fl\left(X\#^n\overline{\mathbb{C}P^2},U\right)\ne 0$ and $\fl\left(\#^n\overline{\mathbb{C}P^2}\right)\ne 0$.
	\end{proposition}

	\begin{remark}
		Ren and Willis conjectured in \cite[Conjecture 6.19]{RW26} that every Stein $2$-handlebody has nonvanishing Khovanov skein lasagna modules. The above proposition partially confirms this conjecture on the link Floer side, provided that the boundary link is a  local unknot rather than the empty link.
	\end{remark}
	
	Another family of nonvanishing results comes from the knot trace $X_n(K)$,  which  is constructed by attaching an $n$-framed $2$-handle to $B^4$ along $K\subset S^3$.

	\begin{proposition}\label{thm1.0}
		Let \(K\) be a knot in \(S^3\), and set $\widehat{\nu}(K)=|\varepsilon(K)|\bigl(2\nu(K)-1\bigr)$,
		where \(\varepsilon(K)\in\{-1,0,1\}\) and \(\nu(K)\in\mathbb{Z}\) are concordance invariants of \(K\) from \cite{Hom14,OS11}. Then, for every integer \(n\le \widehat{\nu}(K)+0^{\widehat{\nu}(K)}-1\), we have
		\[
		\fl(X_n(K),U)\ne 0,
		\]
		where \(U\subset\partial X_n(K)\) is a local unknot, and \(0^0=1\), \(0^k=0\) for \(k\ne0\).
	\end{proposition}
	
	\begin{remark}
		Let $\m{TB}(K)$ be the maximal Thurston-Bennequin number, then for any $n\le \m{TB}(K)-1$, we have $X_n(K)$ is a Stein domain. Hence from Proposition \ref{thmstein} we have $\fl(X_n(K),U)\ne 0$. However, \cite[Theorem 1]{Pla04} showed that $\m{TB}(K)\le 2\tau(K)-1$, hence  $\m{TB}(K)-1\le \widehat{\nu}(K)+0^{\widehat{\nu}(K)}-1$. Thus Proposition \ref{thm1.0}  gives more nonvanishing results for knot traces than   Proposition \ref{thmstein}.
	\end{remark}

	We can generalize the definition of the  canonical Floer lasagna filling mentioned above, to any $4$-manifold $X$ with empty boundary link. We still denote it as $\mathcal{F}_X$. This canonical Floer lasagna filling is functorial under the inclusion map. If $i:X\to X'$ is an inclusion between $4$-manifolds, then $i$ induces a map $i_*$ on the Floer lasagna module such that $i_*\left([\mathcal{F}_X]\right)=\left[\mathcal{F}_{X'}\right]$. Thus $\left[\mathcal{F}_{X'}\right]\ne 0$ implies $\left[\mathcal{F}_{X}\right]\ne 0$ and $\fl(X)\ne 0$.
	
	Therefore, since \(\left[\mathcal{F}_{\#^n\overline{\mathbb{C}P^2}}\right]\ne 0\) by Proposition \ref{thmstein}, it follows that for every \(4\)-manifold \(X\) that embeds into some \(\#^n\overline{\mathbb{C}P^2}\), the Floer lasagna module \(\fl(X)\) is nonzero. In particular, using the folklore trace embedding lemma, we obtain the following.
	
	\begin{proposition}\label{thmknottrace}
		If $K$ is a slice knot in $S^3$, then for any integer $n\le 0$, we have $
		\fl(X_n(K))\ne 0$. 
	\end{proposition}
	
	\begin{remark}\label{31knot}
		As observed in \cite[Example 6.2]{RW26}, $X_n(3_1)$ can be embedded in some $\#^k\overline{\mathbb{C}P^2}$ for any $n\le -4$, where $3_1$ denotes the left-handed trefoil. Hence we have $\fl\left(X_n(3_1)\right)\ne 0$ for any $n\le -4$.
	\end{remark}

The boundary of \(X_n(K)\) is the surgery manifold \(S^3_n(K)\), and the core of the surgery torus in \(S^3_n(K)\) is the dual knot \(K_n\). The remaining results in this paper mainly concern the Floer lasagna module of \((X_n(K),K_n)\).

Equivalently, \(K_n\) can be viewed as the meridian of \(K\), which bounds the cocore disk of the \(2\)-handle. Consequently, \(K_n\) bounds a cocore disk \(D^{cc}\) in \(X_n(K)\). Consider the capped-off Seifert surface \(\widehat{F}\) of \(K\), namely, the Seifert surface of \(K\) together with the core disk. It is easy to see that \(\mathrm{H}_2(X_n(K),K_n)\) is generated by \([D^{cc}]\) and \([\widehat{F}]\).
	
	As observed by Chen in \cite{chen2022floer}, the Floer lasagna module decomposes according to the relative homology classes in \(\text{H}_2(X_n(K),K_n)\), namely
	\[
	\fl(X_n(K),K_n)
	=
	\bigoplus_{a,b\in\mathbb Z}
	\fl\left(X_n(K),K_n,a[D^{cc}]+b[\widehat{F}]\right).
	\]
For any integer \(k\), we set \(\alpha_k=[D^{cc}]+k[\widehat{F}]\) for simplicity, and primarily focus on the nonvanishing of \(\fl(X_n(K),K_n,\alpha_k)\).

The proof of Proposition \ref{thm1.0} relies on a nonvanishing result for Heegaard Floer cobordism maps, derived from the mapping cone formula in \cite{OS08}. Analogously, by using the knot Floer mapping cone formula established by Hedden and Levine in \cite{HL24}, we obtain a nonvanishing result for the knot Floer cobordism maps associated with the dual knot. This allows us to extend Proposition \ref{thm1.0} to the dual knot setting as follows.
	
	\begin{theorem}\label{thm1.2}
		Let $K$ be a knot in $S^3$, then for any nonzero integer $n$,$$\dim \mathcal{FL}(X_n(K),K_n,\alpha_0)\ge \max\{2\tau(K)-n,0\}.$$In particular, if $n\le 2\tau(K)-1$ and $n\ne 0$, then  $\mathcal{FL}(X_n(K),K_n,\alpha_0)\ne 0$.
	\end{theorem}

	For other $k$, we use some special case of surgery exact triangle to obtain the following results.

	\begin{theorem}\label{thm1.1}
		Let $K$ be a knot in $S^3$, then we have 
		\begin{enumerate}
			\item for any integer $n$,  $$\dim \mathcal{FL}(X_n(K),K_n,\alpha_{-1})\ge \dfrac{\dim \widehat{\hfk}(K)-1}{2}.$$

			\item for any integer $n$, $$\dim \mathcal{FL}(X_n(K),K_n,\alpha_{-2})\ge \dfrac{\dim \widehat{\m{HFK}}(K_{2,2n+1})-1}{2}.$$
		\end{enumerate}
	\end{theorem}

	For  L-space knots, we obtain the following stronger dimension bounds.

	\begin{theorem}\label{thm1.5}
		Let $K\subset S^3$ be an   L-space knot, with Seifert genus $g(K)$, then for any integer $n\ge \max\{2g(K),1\}$ and  $k\ge 2$, we have $$\dim \fl(X_n(K),K_n,\alpha_{-k})\ge n.$$Hence $\fl(X_n(K),K_n)$ is infinite dimensional.
	\end{theorem}

	As a corollary of Theorem \ref{thm1.2} and \ref{thm1.1}, we have the following,
	
	\begin{corollary}
		For any knot \(K\subset S^3\) and any integer \(n\),  $\mathcal{FL}(X_n(K),K_n)=0$ if and only if $K$ is the unknot and $n=0$.
	\end{corollary}

	\begin{proof}
		If  $\mathcal{FL}(X_n(K),K_n)=0$,  then  $\mathcal{FL}(X_n(K),K_n,\alpha_k)=0$ for every $k$. From Theorem \ref{thm1.1} (1) we know $\dim \widehat{\hfk}(K)=1$ hence $K$ is the unknot. And similarly from Theorem \ref{thm1.1} (2) we have $K_{2,2n+1}=U_{2,2n+1}$ is also an unknot, then $n=-1$ or $0$. Now apply Theorem \ref{thm1.2} and note that $2\tau(U)-1=-1$, hence if $n=-1$ then $\fl(X_n(K),K_n,\alpha_0)\ne 0$ which is a contradiction, then we know $n$ must be zero.
		
		On the other hand, if $K=U$ and $n=0$, then it is easy to see $U_0$ intersects a capping disk of $U$ geometrically once, hence $\fl(X_0(U),U_0)=0$ follows from \cite[Proposition 4.1]{chen2022floer}.
	\end{proof}
	
	\begin{remark}
		In \cite{RW26}, they proved that the skein lasagna modules have a vanishing result as follows. Let $K\subset S^3$ be a knot, then for any $n\ge -\m{TB}(m(K))$, the skein lasagna module $\mathcal{S}^2_0(X_n(K),L)=0$ for any $L\subset \partial X_n(K)$. In contrast, our result shows that similar vanishing result will never happen in Floer lasagna module when $K$ is nontrivial or $n$ is nonzero, which  reflects the differences between the Floer lasagna module and the skein lasagna module.
	\end{remark}

	\noindent	\textbf{Organization of the paper.} In Section \ref{sec2} we provide some background on link Floer homology. In Section  \ref{sec3} we give the definition of Floer lasagna modules and use some basic properties of Floer lasagna module and Heegaard Floer homology to first prove Theorem \ref{thmgeneralW}, then apply it to prove Proposition  \ref{thmstein}, \ref{thm1.0} and \ref{thmknottrace}. In Section \ref{sec4}, we first use the dual knot mapping cone formula to calculate the rank of a specific  annulus cobordism map, which may be of independent interest, then use the calculation to prove Theorem \ref{thm1.2}. In Section \ref{sec5}, we construct a special surgery triad to calculate more annulus cobordism maps and prove Theorem \ref{thm1.1}, then restrict to the  L-space knot case to prove Theorem \ref{thm1.5}.

	\noindent	\textbf{Acknowledgements.} The author would like to thank his advisor Fan Ye for his continuous guidance and assistance throughout this project. The author would also like to thank  Daren Chen and Cai Zeng for sharing their thoughts, reading an earlier draft of this paper and providing valuable comments. In addition, the author would like to thank Yuan Xu for his advice on writing conventions.
	
	\noindent	\textbf{AI Declaration.} The computation of the dimension of knot Floer homology for L-space knots in Section \ref{calculation L space} was obtained through discussions between the author and ChatGPT 6 Astra. ChatGPT 6 Astra also assisted with polishing the manuscript, including drawing figures, checking grammar and improving the English. No AI tools were used for any other aspects of this work. The author takes full responsibility for the	mathematical arguments and final text.

	\section{Link Floer homology}\label{sec2}
	
	In this section, we provide some background on  link Floer homology.
	
	\begin{definition}
		A  \textit{multi-based link} in a  $3$-manifold $Y$ is an ordered triple $\mathbb{L}=(L,\textbf{w},\textbf{z})$ where $L$ is  an oriented
		link,   $\textbf{w}$ and $\textbf{z}$ are two disjoint collections of basepoints on $L$, such that each component of $L$ contains at least two basepoints, and the basepoints in $\textbf{w}$ and $\textbf{z}$ alternate along each component of $L$. Furthermore, each component of $Y$ contains at least one component of $L$.
	\end{definition} 
	
	In this paper, we primarily consider doubly based knots. We suppress the basepoints from the notation when no confusion can arise. Hence we often  write $L$ rather than $\mathbb{L}$ or $(L,\textbf{w},\textbf{z})$ for simplicity, if each component of $\mathbb{L}$ has exactly two basepoints.

	For any pair $(Y,\mathbb{L})$, we can associate a link Floer homology $\widehat{\hfl}(Y,\mathbb{L})$ over the field $\mathbb{F}_2$ defined by Zemke in \cite{Zem19a}. And we set $\widehat{\hfl}(\varnothing,\varnothing)=\mathbb{F}_2$. If each component of $\mathbb{L}$ has exactly two basepoints, then $\widehat{\hfl}(Y,\mathbb{L})$ agrees with the link Floer homology group $\widehat{\hfl}(Y,L)$ defined by Ozsv\'ath and Szab\'o in \cite{OS08b}.
	
	When \(Y=S^3\) and each component of $\mathbb{L}$ carries exactly two basepoints,  the link Floer homology $\widehat{\hfl}(L)=\widehat{\hfl}(S^3,\mathbb{L})$ admits a $\mathbb{Z}$-grading called the Maslov grading $M$ and a $\dfrac{1}{2}\mathbb{Z}^{|L|}$-grading called the Alexander multigrading $\textbf{A}=(A_1,\cdots,A_{|L|})$, and a decomposition associated to these gradings: $$\widehat{\hfl}(L)=\bigoplus_{\substack {M\in \mathbb{Z}\\ \textbf{A}\in \frac{1}{2}\mathbb{Z}^{|L|}}}\widehat{\hfl}_M(L,\textbf{A}),$$we often collapse the Alexander multigrading to a single Alexander grading $A$ as $$\widehat{\hfl}_M(L,A)=\bigoplus_{|\textbf{A}|=A}\widehat{\hfl}_M(L,\textbf{A}),$$where $|\textbf{A}|=A_1+\cdots +A_{|L|}$.

	For general $3$-manifold and rationally null-homologous multi-based link, the link Floer homology has gradings $\m{gr}_{\textbf{w}}$ and $\m{gr}_{\textbf{z}}$ defined in \cite{Zem19b}, and the Alexander grading $A$ is determined by $\dfrac{1}{2}\left(\m{gr}_{\textbf{w}}-\m{gr}_{\textbf{z}}\right)$, which agrees with the collapsed Alexander grading mentioned above when $Y=S^3$. And in that case we can view $\m{gr}_{\textbf{w}}$ as the Maslov grading $M$ with some grading shift, one can see \cite[Section 1.5]{Zem19b} for more details.
	
	If $\mathbb{K}=(K,w,z)$ in $Y$ is a knot with two basepoints, we  also simply denote $\widehat{\hfk}(Y,K)$ as $\widehat{\hfl}(Y,\mathbb{K})$, which is not an abuse of notation since $\widehat{\hfl}(Y,\mathbb{K})$ is indeed isomorphic to the Ozsv\'ath-Szab\'o's knot Floer homology.  We will talk more about the Alexander grading of this knot Floer homology group in Section \ref{mappingconeformula}.
	
	To define link Floer cobordism maps, we equip cobordism surfaces with decorations,  which was first introduced by Juh\'asz in \cite{Juh16} for hat version of link Floer homology and then generalized to minus and infinity version by Zemke in \cite{Zem19a}. 
	
	\begin{definition}
		We say a pair \((W, \mathcal{F})\) is a \textit{decorated link cobordism} between two 3-manifolds with multi-based links, and write \((W, \mathcal{F}) \colon (Y_1, L_1, \mathbf{w}_1, \mathbf{z}_1) \to (Y_2, L_2, \mathbf{w}_2, \mathbf{z}_2)\), if the following are satisfied.
		\begin{enumerate}
			\item \(W\) is a $4$-dimensional cobordism from \(Y_1\) to \(Y_2\).
			\item \(\mathcal{F} = (\Sigma, \mathcal{A})\), where \(\Sigma\) is an oriented surface, properly embedded in \(W\), with \(\partial \Sigma = -L_1 \sqcup L_2\).
			\item \(\mathcal{A} \subseteq \Sigma\) is a properly embedded 1-manifold dividing \(\Sigma\) into two  subsurfaces with disjoint interiors \(\Sigma_{\mathbf{w}}\) and \(\Sigma_{\mathbf{z}}\), which meet along \(\mathcal{A}\).
			\item Each component of \(L_i \setminus \mathcal{A}\) contains exactly one basepoint.
			\item \(\mathbf{w}_1 \cup \mathbf{w}_2 \subseteq \Sigma_{\mathbf{w}}\) and \(\mathbf{z}_1 \cup \mathbf{z}_2 \subseteq \Sigma_{\mathbf{z}}\).
		\end{enumerate}
	\end{definition}
	
	\begin{example}
		We  mainly use the following  two types of decorated surfaces in this paper, 
		\begin{enumerate}
			\item If $\Sigma=K\times [0,1]$ is a product annulus, and $K$ has basepoints $w,z$. Choose points \(p,q\in K\setminus\{w,z\}\), one in each component of \(K\setminus\{w,z\}\). The product arcs \(\mathcal A=\{p,q\}\times[0,1]\)  are called   \textit{canonical divides on the annulus}.
			
			\item If $D$ is a disk with boundary of $K$, take $\mathcal{A}_D$ a single arc seprating $D$ into $w$ and $z$ two regions as the  dividing, it is called  \textit{canonical divide on the disk}.
		\end{enumerate}
		\begin{figure}[htbp]
			\centering
			\includegraphics[width=0.5\linewidth]{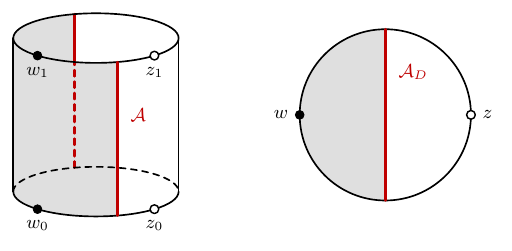}
			\caption{Left: canonical divides on annulus. Right: canonical divide on disk. }
			\label{fig:canonicaldivides}
		\end{figure}
		Throughout the rest of the paper, unless otherwise specified, we always take canonical divides on the annulus and the disk. 
	\end{example}

	A decorated link cobordism $(W,(\Sigma,\mathcal{A})):(Y_1,\mathbb{L}_1)\to (Y_2,\mathbb{L}_2)$ induces,  for each $\spc$ structure $\mathfrak{s}$ on $W$, a link Floer cobordism map \begin{equation}\label{linkfloercob}
		\widehat{F}_{W,(\Sigma,\mathcal{A}),\mathfrak{s}}: \widehat{\hfl}(Y_1,\mathbb{L}_1)\to \widehat{\hfl}(Y_2,\mathbb{L}_2).
	\end{equation}

	When $c_1(\mathfrak{s}|_{Y_i})$ is torsion, and each component of $L_i$ is rationally null-homologous,  we can describe the Alexander grading shift of the link Floer cobordism map in \eqref{linkfloercob}. Following  \cite[Section 3.2]{HHSZ26}, suppose that $L_1\subset Y_1$ and $L_2\subset Y_2$ are two $\ell$-component links. Write $L_j^1,\cdots,L_j^\ell$ for the components of $L_j$. Suppose that each component is rationally null-homologous with two basepoints and further that we have choices $T_1=(T_1^i)_{1\le i\le \ell}$ and $T_2=(T_2^i)_{1\le i\le \ell} $ of collections of rational $2$-chains,  such that \(\partial T_j^i = -L_j^i\), namely,  \(T_j^i\) is a \textit{rational Seifert surface} for \(L_j^i\).
	
	Suppose that \((W, \Sigma) : (Y_1, L_1) \to (Y_2, L_2)\) is an oriented link cobordism. If \(\Sigma\) has components \(\Sigma_1, \cdots, \Sigma_\ell\), such that \(\partial \Sigma_i = -L_1^i \cup L_2^i\), then  the choices of rational Seifert surfaces give a unique lift of \([\Sigma_i] \in \m{H}_2(W, \partial W, \mathbb{Q})\) to \(\m{H}_2(W, \mathbb{Q})\). Write \([\widehat{\Sigma}_i]\) for this class, then \([\widehat{\Sigma}]\) decomposes as a sum \([\widehat{\Sigma}_1] + \dots + [\widehat{\Sigma}_\ell]\).
	
	For our use, we restrict to the case where $\Sigma$ is a collection of annuli, and we pick the canonical divides on each annulus and these divides together form a decoration on $\Sigma$. Now we can state the grading shift result as follows.
	
	\begin{proposition}[{\cite[Proposition 3.1]{HHSZ26}}]\label{alexandershift}\quad
		
		\begin{enumerate}
			\item Suppose that \(L \subset Y\) is an \(\ell\)-component link, each of whose components is rationally null-homologous. Let \(T = (T^i)_{1 \le i \le \ell}\) denote a collection of rational Seifert surfaces for the components of \(L\). Then there is a well-defined \(\ell\)-component Alexander grading \(A_T = (A_{T, 1}, \dots, A_{T, \ell})\) on \(\widehat{\hfl}(Y, L, \mathfrak{s})\). If \(c_1(\mathfrak{s})\) is torsion, the grading is independent of the choice of \(T\).
			\item Suppose \((W, \Sigma) : (Y_1, L_1) \to (Y_2, L_2)\) is a link cobordism, consisting of \(\ell\) disjoint annuli, and let \(\mathcal{F} = (\Sigma, \mathcal{A})\) described as above. Suppose that \(T_1\) and \(T_2\) are collections of rational Seifert surfaces, as above. Then
			\[
			A_{T_2, i}(\widehat{F}_{W, \mathcal{F}, \mathfrak{s}}(\mathbf{x})) - A_{T_1, i}(\mathbf{x}) = \frac{\langle c_1(\mathfrak{s}), [\widehat{\Sigma}_i] \rangle - [\Sigma] \cdot [\widehat{\Sigma}_i]}{2}.
			\]
		\end{enumerate}
	\end{proposition}

	\begin{remark}\label{emptyremark}
		As implicitly used in the proof of \cite[Lemma 3.2]{HHSZ26}, similar result still holds when $1\le |L_2|<\ell$. In that case, we need $\Sigma$ to   consist of $|L_2|$-annuli and $\ell-|L_2|$ disks. More precisely, we can view $L_2$ as $L_2^1\cup \cdots\cup L_2^{k}\cup \varnothing\cup\cdots\cup \varnothing$, where $k=|L_2|$ and we have $\ell-k$ empty links. And the Alexander grading of the empty knot is set to be zero.
	\end{remark}

	\section{Floer lasagna modules}\label{sec3}

	\begin{definition}
		Let \(W\) be a smooth oriented 4-manifold, and \(\mathbb{L}=(L,\mathbf{w},\mathbf{z})\) a multi-based link in \(\partial W\). A \textit{Floer lasagna filling}  of \((W,\mathbb{L})\) is a tuple $$\mathcal{F}=(\{B_i,\mathbb{L}_i,v_i\}_{i\in I},(\Sigma,\mathcal{A}))$$for some nonempty finite indexing set \(I\), such that for all $i\in I$:
		\begin{enumerate}
			\item   The \(B_i\) are pairwise disjoint smoothly embedded \(4\)-balls in \(W^\circ\), called input balls. 
			
			\item   \(\mathbb{L}_i=(L_i,\mathbf{w}_i,\mathbf{z}_i)\) is a multi-based link in $\partial B_i$.
			
			\item A label \(v_i\in\widehat{\hfl}(\partial B_i,\mathbb{L}_i)\) for each \(i\), where the link Floer homology groups are taken with coefficients \(\mathbb{F}_2\).

			\item  \((\Sigma,\mathcal{A})\) is a decorated oriented properly embedded surface in \(W\backslash(\sqcup_i B_i^\circ)\) such that  \(\Sigma\cap\partial W=L\) and \(\Sigma\cap\partial B_i=L_i\), whose induced orientations agree with those of \(L\) and \(-L_i\).
		\end{enumerate}
	\end{definition}
	
	\begin{example}
		As mentioned in Section \ref{sec1},	for any pair $(W,\varnothing)$, we have a \textit{canonical Floer lasagna filling} \begin{equation}\label{canonicalfilling}
			\mathcal{F}_W=\left(\{B,\mathbb{U},1\},(D,\mathcal{A}_D)\right),
		\end{equation}where $B$ is a small input ball in the interior of $W$, and $\mathbb{U}=(U,w,z)\subset \partial B$ is an unknot with two basepoints, $1\in \widehat{\hfk}(S^3,\mathbb{U})\cong \mathbb{F}_2$. And $D$ is a small disk bounded by $U$ in $W\setminus B^\circ$, and $\mathcal{A}_D$ is the canonical divides on disk $D$.
		
		Similarly, if $W$ has nonempty boundary $Y$ and $U\subset Y$ is a local unknot, we also have a \textit{canonical Floer lasagna filling} \begin{equation}\label{canonicalfillingU}
			\mathcal{F}_{W,U}=\left(\{B,\mathbb{U},1\}, (U\times I,\mathcal{A})\right),
		\end{equation}which is constructed as follows. Since $Y$ is nonempty, we can assume $W$ admits a handle decomposition with a single $0$-handle and without $4$-handle. The $4$-ball $B$ is the single $0$-handle.  Now take a local unknot $U$ away from the attaching spheres of the handles, then track its image we have a product annulus cobordism  $\left(W\setminus B^\circ, (U\times I,\mathcal{A})\right):(S^3,U)\to (Y,U)$, where $\mathcal{A}$  is the canonical divides on the annulus.
	\end{example}
	
	\begin{definition}\label{simprime}
		Let \(\mathcal{F} = \left( \left\{B_i, \mathbb{L}_i, v_i\right\}_{i \in I},(\Sigma, \mathcal{A})\right)\) and \(\mathcal{F}' = \left( \left\{B'_j, \mathbb{L}'_j, v'_j\right\}_{j \in J},(\Sigma', \mathcal{A}')\right)\) be lasagna fillings of \((W, \mathbb{L})\). We say \(\mathcal{F} \sim' \mathcal{F}'\) if the following hold:
		\begin{enumerate}
			\item[(1)] we can partition \(J = \sqcup_{i \in I} J_i\) such that for each  \(i \in I\) we have \(\sqcup_{j \in J_i} B'_j \subseteq B_i\);
			\item[(2)] \((\Sigma', \mathcal{A}')\) is isotopic, relative to the boundary, to \((\Sigma, \mathcal{A})\) in \(W \setminus (\sqcup_{i \in I} B_i^\circ)\);
			\item[(3)] for each \(i \in I\), \(\Sigma' \cap \partial B_i\) is isotopic to \(\mathbb{L}_i\) where the basepoints are determined by the subsurfaces \(\Sigma_{\mathbf{w}}'\) and \(\Sigma_{\mathbf{z}}'\);
			\item[(4)] for each \(i \in I\), the map \(\widehat{F}_{\mathcal{F}_i}\) coming from the link cobordism
		\begin{equation}\label{Fi}
			\mathcal{F}_i = \left(B_i\setminus \left(\sqcup_{j \in J_i} B'_j\right)^\circ, (\Sigma', \mathcal{A}') \cap B_i  \right) : \left(\sqcup_{j \in J_i} \partial B_j',\sqcup_{j \in J_i}\mathbb{L}_j'\right) \to \left(\partial B_i,\mathbb{L}_i\right)
		\end{equation}
			has the property \(\widehat{F}_{\mathcal{F}_i}(\otimes_{j\in J_i} v_j') = v_i\).
		\end{enumerate}
		We then define an equivalence relation \(\sim\) on \(\mathbb{F}\{\text{lasagna fillings of } (W, \mathbb{L})\}\)  generated by the local replacement relations $\sim'$ and multilinearity in the \(v_i\). That is, if \(\mathcal{F} = ( \{B_i, \mathbb{L}_i, v_i\}_{i \in I},(\Sigma, \mathcal{A})) \) is a lasagna filling and for some \(j \in I\) we have \(v_j = w_j + z_j\), then
		\[
		\mathcal{F} \sim ( \{B_i, \mathbb{L}_i, w_i\}_{i \in I},(\Sigma, \mathcal{A})) + ( \{B_i, \mathbb{L}_i, z_i\}_{i \in I},(\Sigma, \mathcal{A}))
		\]
		where \(w_i = z_i = v_i\) for \(i \neq j\).

	\end{definition}
	
	Using this equivalence relation we can define the Floer lasagna module of the pair $(W,\mathbb{L})$
	
\begin{definition}
	We define the \textit{Floer lasagna module} of $(W,\mathbb{L})$ to be 
	\begin{equation}
		\mathcal{FL}(W,\mathbb{L}):=\mathbb{F}\left\{\text{Floer lasagna fillings of } (W,\mathbb{L})\right\}/\sim.
	\end{equation}
	We use \([\mathcal{F}]\) to denote the class of the filling \(\mathcal{F}\) in the Floer lasagna module.
\end{definition}

	As observed by Chen in \cite{chen2022floer}, we can decompose $\mathcal{FL}(W,\mathbb{L})$ by the relative homology class in $\text{H}_2(W,L)$, hence we have $$\mathcal{FL}(W,\mathbb{L})=\bigoplus_{\alpha\in\text{H}_2(W,L)}\mathcal{FL}(W,\mathbb{L},\alpha).$$
	
	Chen also observed an evaluation map from the Floer lasagna module to the link Floer homology group. We can further refine Chen's evaluation map by incorporating $\spc$ structures.
	
	\begin{proposition}
		For any $\m{Spin}^c$ structure  $\mathfrak{u}$ on $W$,  there is an evaluation map \begin{equation*}
			\mathrm{ev}_{\mathfrak{u}}: \mathcal{FL}(W,\mathbb{L})\to \widehat{\hfl}(\partial W,\mathbb{L}).
		\end{equation*}For any Floer lasagna filling $[\mathcal{F}]=\left[(\{B_i,\mathbb{L}_i,v_i\}_{i\in I},(\Sigma,\mathcal{A}))\right]$, it is defined by \begin{equation}\label{ev}
			\mathrm{ev}_{\mathfrak{u}}\left([\mathcal{F}]\right):=\widehat{F}_{W\setminus \sqcup_i B_i^\circ,(\Sigma,\mathcal{A}),\mathfrak{u}}(\otimes_i v_i),
		\end{equation}where $\widehat{F}_{W\setminus \sqcup_i B_i^\circ,(\Sigma,\mathcal{A}),\mathfrak{u}}: \widehat{\hfl}\left(\sqcup \partial B_i,\sqcup_i\mathbb{L}_i\right)\cong \otimes_i\widehat{\hfl}(\partial B_i,\mathbb{L}_i)\to\widehat{\hfl}(\partial W,\mathbb{L}) $ is the link Floer cobordism map induced from $\left(W\setminus \sqcup_i B_i^\circ ,(\Sigma,\mathcal{A}),\mathfrak{u}\right)$. We denote $\m{ev}=\displaystyle \sum_{\mathfrak{u}\in \m{Spin}^c(W)}\m{ev}_{\mathfrak{u}}$.
	\end{proposition}
	
	\begin{remark}
	Here we use the fact that $\m{Spin}^c$ structures extend uniquely over \(4\)-balls. Hence, we can canonically identify $\m{Spin}^c(W)$ with $\m{Spin}^c(W\setminus \sqcup_i B_i^\circ)$, and we view $\mathfrak{u}$ as an element of the latter.
	\end{remark}

	\begin{proof}
		We  show that this evaluation map is well-defined. The linearity of the cobordism maps ensures compatibility with the multilinearity relations. Therefore, it suffices to show that $\mathrm{ev}_{\mathfrak{u}}([\mathcal{F}])= \mathrm{ev}_{\mathfrak{u}}([\mathcal{F}'])$ whenever $\mathcal{F}\sim'\mathcal{F}'$. 	We adopt the notation from Definition \ref{simprime}. Then we have 
		\begin{align*}
			\m{ev}_{\mathfrak{u}}\left([\mathcal{F}]\right)&=\widehat{F}_{W\setminus \sqcup_i B_i^\circ,(\Sigma,\mathcal{A}),\mathfrak{u}}(\otimes_{i\in I} v_i),\\
			\m{ev}_{\mathfrak{u}}\left([\mathcal{F}']\right)&=\widehat{F}_{W\setminus \sqcup_j (B_j')^\circ,(\Sigma',\mathcal{A}'),\mathfrak{u}}(\otimes_{i\in J} v_j').
		\end{align*} 
	
		Now note that 
		\begin{equation}\label{decomp}
			W\setminus \sqcup_j (B_j')^\circ=\left(W\setminus \sqcup_i B_i^\circ\right)\sqcup_{i\in I}\left(B_i\setminus\sqcup_{j\in J_i}(B_j')^\circ\right).
		\end{equation} 
		Furthermore, $\spc\left(W\setminus \sqcup_j (B_j')^\circ\right)\cong \spc\left(W\setminus \sqcup_i B_i^\circ\right)$. Indeed, since $\m{H}^2\left(B_i\setminus\sqcup_{j\in J_i}(B_j')^\circ\right)=0$, any $\spc$ structure on the boundary uniquely extends over the interior pieces $\sqcup_{i \in I} \left(B_i\setminus \sqcup_{j\in J_i}(B_j')^\circ\right)$. 
		
		We recall the composition law from \cite[Theorem B]{Zem19a}: if $(W_1,(\Sigma_1,\mathcal{A}_1))$ is followed by $(W_2,(\Sigma_2,\mathcal{A}_2))$, and  $\mathfrak{s}_1\in\spc(W_1)$ and $\mathfrak{s}_2\in\spc(W_2)$ have matching restrictions to the common boundary, then 
		\begin{equation}\label{linkfloercompo}
			\sum_{\substack{\mathfrak{t}\in\spc(W_1\cup W_2) \\  \mathfrak{t}|_{W_i}=\mathfrak{s}_i}} \widehat{F}_{W_1\cup W_2,(\Sigma_1,\mathcal{A}_1)\cup (\Sigma_2,\mathcal{A}_2),\mathfrak{t}}= \widehat{F}_{W_1,(\Sigma_1,\mathcal{A}_1),\mathfrak{s}_1}\circ \widehat{F}_{W_2,(\Sigma_2,\mathcal{A}_2),\mathfrak{s}_2}.
		\end{equation}
		Because of the uniqueness of the extension of the \(\spc\) structure mentioned above, the sum on the left-hand side of \eqref{linkfloercompo} consists of exactly one term when applied to our decomposition \eqref{decomp}. Therefore, we obtain  
		\begin{align*}
			\widehat{F}_{W\setminus \sqcup_j (B_j')^\circ,(\Sigma',\mathcal{A}'),\mathfrak{u}}(\otimes_j v_j')&=\widehat{F}_{W\setminus \sqcup_i B_i^\circ,(\Sigma,\mathcal{A}),\mathfrak{u}}\circ \widehat{F}_{\sqcup_{i\in I}\mathcal{F}_i}\left(\otimes_{i\in I}\left(\otimes_{j\in J_i}v_j'\right)\right)\\
			&=\widehat{F}_{W\setminus \sqcup_i B_i^\circ,(\Sigma,\mathcal{A}),\mathfrak{u}}\left(\otimes_{i\in I}v_i\right),
		\end{align*}
	where the second equality uses \(\widehat{F}_{\mathcal{F}_i}(\otimes_{j\in J_i} v_j') = v_i\) from Definition \ref{simprime} (and $\mathcal{F}_i$ is defined as in \eqref{Fi}). This yields $\mathrm{ev}_{\mathfrak{u}}([\mathcal{F}])= \mathrm{ev}_{\mathfrak{u}}([\mathcal{F}'])$, completing the proof.
	\end{proof}
	
	\begin{remark}
		In link Floer homology we set $\widehat{\hfl}(\varnothing,\varnothing)\cong \mathbb{F}_2$, but $\widehat{\hfl}(Y,\varnothing)$ is not defined. Hence when the boundary link $\mathbb{L}$ is empty, we can only define the evaluation map when $W$ is closed.
	\end{remark}
	
	The strategy to prove our nonvanishing result  is to first find a special Floer lasagna filling representative, then calculate the rank of the link Floer cobordism map it induces, and finally from \eqref{ev} to deduce \begin{equation}\label{evin}
		\dim \mathcal{FL}(W,\mathbb{L})\ge \text{rank} \text{ ev}\ge \text{rank } \widehat{F}_{(\Sigma,\mathcal{A})}.
	\end{equation}
	
We can also define inclusion maps and gluing maps of Floer lasagna modules analogously to skein lasagna modules in \cite{RW26}.

\begin{definition}[Inclusion map]
	An embedding \(i:X\hookrightarrow (X')^\circ\) induces a natural map $i_*$ between their Floer lasagna modules,
	\begin{align*}
		i_*:\fl(X)\to \fl(X'),\quad [\mathcal{F}]\mapsto [\mathcal{F}],
	\end{align*}
	since each Floer lasagna filling of $(X,\varnothing)$ can be naturally viewed as a Floer lasagna filling of $(X',\varnothing)$ through the inclusion.
	
	It is easy to see that the canonical Floer lasagna filling \eqref{canonicalfilling} is functorial under the inclusion map, namely,
	\begin{equation}\label{funcinc}
		i_*\left([\mathcal{F}_X]\right)=[\mathcal{F}_{X'}].
	\end{equation}
\end{definition} 

More generally, we can extend the inclusion map to a gluing map as follows.

	\begin{definition}[Gluing map]Let \(Y=\partial W\), and let $(W',(\Sigma',\mathcal{A}')):(Y,\mathbb{L})\to (Y',\mathbb{L}')$ be a decorated link cobordism, we can define a gluing map as follows. \begin{align*}
			\mathsf{G}_{W',(\Sigma',\mathcal{A}')} :\fl(W,\mathbb{L})&\to \fl(W\cup W',\mathbb{L}')\\
			[\mathcal{F}]=\left[(\{B_i,\mathbb{L}_i,v_i\}_{i\in I},(\Sigma,\mathcal{A}))\right]&\mapsto \left[(\{B_i,\mathbb{L}_i,v_i\}_{i\in I}, (\Sigma,\mathcal{A})\cup (\Sigma',\mathcal{A}'))\right],
		\end{align*}it is easy to see this gluing map is well-defined. 
	\end{definition}

	From the composition law of link Floer homology \eqref{linkfloercompo}, we can easily obtain the following,
	
	\begin{proposition}
		Given $\spc$ structures \(\mathfrak{s}\) on \(W\) and \(\mathfrak{s}'\) on \(W'\), \begin{equation}\label{gluingspinc}
			\sum_{\substack{\mathfrak{t}\in\spc(W\cup W') \\ \mathfrak{t}|_W=\mathfrak{s}, \mathfrak{t}|_{W'}=\mathfrak{s}'}} \m{ev}_{\mathfrak{t}} \left(	\mathsf{G}_{W',(\Sigma',\mathcal{A}')} \left([\mathcal{F}]\right)\right)= \widehat{F}_{W',(\Sigma',\mathcal{A}'),\mathfrak{s}'}\left(\m{ev}_{\mathfrak{s}}\left([\mathcal{F}]\right)\right).
		\end{equation}
	\end{proposition}

	\subsection{Proof of Theorem \ref{thmgeneralW}} As a warm-up, we can now give a proof of Theorem \ref{thmgeneralW} using \eqref{evin}.

	Let $W$ be a smooth  $4$-manifold with connected boundary $Y=\partial W$.  We first show that the class of the  canonical Floer lasagna filling $\left[\mathcal{F}_{W,U}\right]$ of $(W,U)$ is nonzero.

	For any  $\mathfrak{s}\in \spc(W)$, under the identification of $\widehat{\m{HF}}(Y)\cong \widehat{\hfk}(Y,U)$, we have the knot Floer cobordism map $\widehat{F}_{W\setminus B^\circ, (U\times I,\mathcal{A}),\mathfrak{s}}$ can be identified with $\widehat{F}_{W\setminus B^\circ ,\mathfrak{s}}$ by comparing the definitions of two cobordism maps from \cite{OS06,Zem19a}.

	Now if there exists an  $\mathfrak{s}\in\m{Spin}^c(W)$ such that the Heegaard Floer cobordism map \begin{equation}\label{HFcobmapnonzero}
		\widehat{F}_{W\setminus B^\circ,\mathfrak{s}}:\widehat{\m{HF}}(S^3)\to \widehat{\m{HF}}(Y)
	\end{equation}is nonzero, then we have $\widehat{F}_{W\setminus B^\circ,( U\times I,\mathcal{A}),\mathfrak{s}}$ is nonzero, in particular, we have \begin{equation}
		\widehat{F}_{W\setminus B^\circ,( U\times I,\mathcal{A}),\mathfrak{s}}(1)\ne 0 \in \widehat{\hfk}(Y,U).
	\end{equation}

	Therefore, by the definition of the evaluation map, the canonical Floer lasagna filling satisfies
	\[
	\m{ev}_{\mathfrak{s}}\left([\mathcal{F}_{W,U}]\right)=\widehat{F}_{W\setminus B^\circ, (U\times I,\mathcal{A}),\mathfrak{s}}(1)\ne 0,
	\]
	which implies $[\mathcal{F}_{W,U}]\ne 0$, and thus $\fl(W,U)\ne 0$.

	Now if $Y=S^3$, consider the $4$-handle map $\left(B^4,(D,\mathcal{A}_D)\right)$ defined in \cite{Zem19a}, which is a cobordism from $(S^3,U)$ to $(\varnothing,\varnothing)$ by capping off $S^3$ by a  $4$-ball and $U$ by a disk $D$. And the induced knot Floer cobordism map $\widehat{F}_{B^4,(D,\mathcal{A}_D)}$ between $\widehat{\hfk}(S^3,U)\cong \mathbb{F}_2$ and $\widehat{\hfl}(\varnothing,\varnothing)\cong \mathbb{F}_2$ is the identity.

	Now as mentioned in Section \ref{sec1}, we can cap off the canonical Floer lasagna filling $\mathcal{F}_{W,U}$ through the $4$-handle map. Then we get a closed $4$-manifold $X$, and the new Floer lasagna filling is nothing but the canonical Floer lasagna filling $\mathcal{F}_{X}$ defined in  \eqref{canonicalfilling}. More precisely,  we have $\mathcal{F}_X$ is obtained from $\mathcal{F}_{W,U}$ by gluing $\left(B^4,(D,\mathcal{A}_D)\right)$, \begin{equation}\label{glu4}
		\mathsf{G}_{B^4,(D,\mathcal{A}_D)}\left([\mathcal{F}_{W,U}]\right)=\left[\mathcal{F}_{X}\right].
	\end{equation}
	
	For the $\spc$ structure $\mathfrak{s}$ in \eqref{HFcobmapnonzero}, we can uniquely extend it to a $\spc$ structure $\overline{\mathfrak{s}}$ on $X=W\cup B^4$. And also from the uniqueness of the extension, for any $\mathfrak{t}\in\spc(X)$ if $\mathfrak{t}|_W=\mathfrak{s}$, then $\mathfrak{t}=\overline{\mathfrak{s}}$. Hence combining \eqref{gluingspinc} and \eqref{glu4}, we obtain \begin{align*}
		\m{ev}_{\overline{\mathfrak{s}}}\left([\mathcal{F}_X]\right)&=\m{ev}_{\overline{\mathfrak{s}}}\left(\mathsf{G}_{B^4,(D,\mathcal{A}_D)}\left([\mathcal{F}_{W,U}]\right)\right)=\widehat{F}_{B^4,(D,\mathcal{A}_D)}\left(\m{ev}_{\mathfrak{s}}\left([\mathcal{F}_{W,U}]\right)\right)\ne 0,
	\end{align*}where we use $\m{ev}_{\mathfrak{s}}\left([\mathcal{F}_{W,U}]\right)\ne 0$ and the $4$-handle map  $\widehat{F}_{B^4,(D,\mathcal{A}_D)}$ is identity.

	Thus we have $\left[\mathcal{F}_{X}\right]\ne 0$, hence $\fl\left(X\right)\ne 0$. Then we finish the proof of Theorem \ref{thmgeneralW}.

	\subsection{Proof of Proposition \ref{thmstein} and Proposition \ref{thmknottrace}} 	We can now combine Theorem \ref{thmgeneralW} and various Heegaard Floer cobordism map results to prove Proposition \ref{thmstein} and Proposition \ref{thmknottrace}.

	If $X$ is a Stein domain with boundary $Y$, then from \cite[Theorem 11.2.2]{GS99}, we know that $X$ admits a special handle decomposition with  one single $0$-handle and several  $1,2$-handles. 
	
	Now take the $4$-ball $B\subset X$ as the $0$-handle, the Stein structure of $X$ induces the standard contact structure $\xi_{\m{std}}$ on $S^3$, and also a  contact structure  $\xi$ on $Y$. Hence we have a Stein cobordism \begin{equation}\label{steincob}
		X\setminus B^\circ :(S^3,\xi_{\m{std}})\to (Y,\xi).
	\end{equation}
	
	Recall for each contact $3$-manifold $(Y,\xi)$ we have an associated contact invariant $c(Y,\xi)\in\widehat{\m{HF}}(-Y)$ and the following functoriality result of Stein cobordism.
	
	\begin{lemma}[{\cite[Theorem 3]{HN22}}]
		If $W:(Y_1,\xi_1)\to (Y_2,\xi_2)$ is a Stein cobordism, then we have \begin{equation}\label{steinfunc}
			\widehat{F}_{{W^\vee},\mathfrak{t}}\left(c(Y_2,\xi_2)\right)=c(Y_1,\xi_1),
		\end{equation}where $W^\vee$ is $W$ turned around and viewed as a cobordism from $-Y_2$ to $-Y_1$, $\mathfrak{t}$  is the canonical $\spc$
		structure of $W$ associated to the Stein structure.
	\end{lemma}
	
	Recall if $(\Sigma,\bm{\alpha},\bm{\beta},z)$ is a Heegaard diagram for $Y$, then $(-\Sigma,\bm{\alpha},\bm{\beta},z)$ is a Heegaard diagram for $-Y$. Hence we can identify the generators of $\widehat{\m{CF}}(Y)$ and $\widehat{\m{CF}}(-Y)$ both as $\mathbb{T}_{\bm{\alpha}}\cap \mathbb{T}_{\bm{\beta}}$ and define a pairing\begin{equation}
		\langle\cdot,\cdot \rangle_Y:\widehat{\m{CF}}(Y)\otimes \widehat{\m{CF}}(-Y)\to \mathbb{F}_2
	\end{equation}through $$\langle \textbf{x},\textbf{y}\rangle_Y:=\begin{cases}
		1 & \text{ if } \textbf{x}=\textbf{y}  \\
		0 & \mathrm{ otherwise}
	\end{cases} .$$This pairing map descends to a bilinear pairing on homology \begin{equation}\label{HFpair}
		\langle\cdot,\cdot \rangle_Y:\widehat{\m{HF}}(Y)\otimes \widehat{\m{HF}}(-Y)\to \mathbb{F}_2.
	\end{equation}And the cobordism maps \begin{align*}
		\widehat{F}_{W,\mathfrak{s}}:&\widehat{\m{HF}}(Y_1)\to \widehat{\m{HF}}(Y_2),\\
		\widehat{F}_{W^\vee,\mathfrak{s}}:&\widehat{\m{HF}}(-Y_2)\to \widehat{\m{HF}}(-Y_1)
	\end{align*}are dual to each other under the pairing \eqref{HFpair}: \begin{equation}\label{dualcob}
		\langle \widehat{F}_{W,\mathfrak{s}}(a), b  \rangle_{Y_2}= \langle a, \widehat{F}_{W^\vee,\mathfrak{s}}(b)  \rangle_{Y_1}, 
	\end{equation}where $a\in \widehat{\m{HF}}(Y_1)$ and $b\in \widehat{\m{HF}}(-Y_2)$, see \cite{OS06,Hed08,Zem21} for more details.

	Now together with \eqref{steincob} and \eqref{steinfunc}, we have $$\widehat{F}_{\left(X\setminus B^\circ\right)^\vee, \mathfrak{t}}\left(c(Y,\xi)\right)=c(S^3,\xi_{\m{std}})=1\in \widehat{\m{HF}}(-S^3)\cong \mathbb{F}_2,$$hence for $1\in \widehat{\m{HF}}(S^3)\cong \mathbb{F}_2$, from \eqref{dualcob}, we have $$\langle \widehat{F}_{X\setminus B^\circ,\mathfrak{t}}(1),c(Y,\xi)\rangle_Y=\langle 1,\widehat{F}_{\left(X\setminus B^\circ\right)^\vee, \mathfrak{t}}\left(c(Y,\xi)\right)\rangle_{S^3}=\langle1,1\rangle_{S^3}=1,$$hence we have $\widehat{F}_{X\setminus B^\circ,\mathfrak{t}}(1)\ne 0$, now apply Theorem \ref{thmgeneralW}, we have $\left[\mathcal{F}_{X,U}\right]\ne 0$ hence $\fl(X,U)\ne 0$, which is a special case of Proposition \ref{thmstein} for $n=0$.

	We now use the following blow-up formula to generalize the nonvanishing result for any other $n$ .
	
	\begin{lemma}[{\cite[Theorem 3.7]{OS06}}]\label{Blowup}Let  $W:Y_0\to Y_1$ be a $4$-manifold cobordism between $3$-manifolds, and $W\#^n\overline{\mathbb{C}P^2}$ be the   blow-up of $W$ at $n$ points. For any $\mathfrak{s}\in\spc(W)$, we can uniquely extend it to a $\spc$ structure $\overline{\mathfrak{s}}$ on $W\#^n\overline{\mathbb{C}P^2}$ such that $\overline{\mathfrak{s}}|_W=\mathfrak{s}$ and $\langle c_1(\overline{\mathfrak{s}}),[E_i]\rangle =1$ for all $1\le i\le n$, where $E_i\subset \overline{\mathbb{C}P^2}$ are the exceptional spheres. Then the two Heegaard Floer cobordism maps are equal, $$\widehat{F}_{W,\mathfrak{s}}=\widehat{F}_{W\#^n\overline{\mathbb{C}P^2},\overline{\mathfrak{s}}}:\widehat{\m{HF}}(Y_0)\to \widehat{\m{HF}}(Y_1).$$
	\end{lemma}

	Now for the canonical $\spc$ structure $\mathfrak{t}$ of Stein domain $X\setminus B^\circ $, we uniquely extend it to a  $\spc$ structure $\overline{\mathfrak{t}}$ on   $(X\setminus B^\circ)\#^n\overline{\mathbb{C}P^2}$ as the lemma above. Hence we have $$\widehat{F}_{(X\setminus B^\circ)\#^n\overline{\mathbb{C}P^2},\overline{\mathfrak{t}}}=\widehat{F}_{X\setminus B^\circ,\mathfrak{t}}\ne 0.$$Thus from Theorem \ref{thmgeneralW} we obtain   $\fl\left(X\#^n \overline{\mathbb{C}P^2},U\right)\ne 0$.

	Take $X=B^4$ as the Stein domain,  then we have $X\#^n \overline{\mathbb{C}P^2}=\#^n \overline{\mathbb{C}P^2}\setminus (B^4)^\circ$, and we can fill it with another $4$-ball to get $\#^n\overline{\mathbb{C}P^2}$. Hence again from Theorem \ref{thmgeneralW}, we have $\left[\mathcal{F}_{\#^n\overline{\mathbb{C}P^2}}\right]\ne 0$ thus $\fl\left(\#^n\overline{\mathbb{C}P^2}\right)\ne 0$. Then we finish the proof of Proposition \ref{thmstein}.

	\begin{remark}
		One may wish to generalize our result to any other Stein domain with boundary $S^3$. However, in \cite{Eli90}, Eliashberg showed that if the boundary of a Stein domain is $S^3$, then it must be diffeomorphic to $B^4$.
	\end{remark}
	
	As mentioned in Section \ref{sec1}, we use the fact that $\left[\mathcal{F}_{\#^n\overline{\mathbb{C}P^2}}\right]\ne 0$ and the trace embedding lemma to prove Proposition \ref{thmknottrace}.
	
	\begin{proof}[Proof of Proposition \ref{thmknottrace}]
		Recall the trace embedding lemma in \cite[Lemma 3.3]{HP19} states that $X_n(K)$ can be smoothly embedded into $W$ if and only if the mirror knot $m(K)$ bounds a smooth slice disk $\Sigma$ in $W\setminus (B^4)^\circ$ with $[\Sigma]\cdot [\Sigma]=n$.

		Now suppose $K$ is a slice knot in $S^3$, then consider the slice disk $D$ bounded by $m(K)$ in $B^4$, for any integer $n\ge 0$, we blow up $n$ distinct interior points on $D$, then we get  $4$-manifold $\#^n \overline{\mathbb{C}P^2}\setminus B^4$. And the proper transform  $\Sigma$ of $D$ is still a disk, now with self intersection number $-n$. Hence from the trace embedding lemma above, we have $X_{-n}(K)$ can be smoothly embedded into $\#^n\overline{\mathbb{C}P^2}$ for any integer $n\ge 0$. We denote this embedding as an inclusion map $i$.
		
		Now from \eqref{funcinc}, the functoriality of canonical filling under inclusion map, together with $\left[\mathcal{F}_{\#^n\overline{\mathbb{C}P^2}}\right]\ne 0$, we have $$i_*\left(\left[\mathcal{F}_{X_{-n}(K)}\right]\right)=\left[\mathcal{F}_{\#^n\overline{\mathbb{C}P^2}}\right]\ne 0,$$hence we have $\left[\mathcal{F}_{X_{-n}(K)}\right]\ne 0$ thus $\fl\left(X_{-n}(K)\right)\ne 0$ for any integer $n\ge 0$ and $K$ slice, which finishes the proof.
	\end{proof}

	\subsection{Proof of Proposition \ref{thm1.0}}We now prove Proposition \ref{thm1.0}, which is an immediate corollary of the following Theorem \ref{HW26cal}. For each $s\in \mathbb{Z}$, there is a  $\mathfrak{v}_s\in\spc(X_n(K))$ such that \begin{equation}
		\left\langle c_1(\mathfrak{v}_s),[\widehat{F}]  \right\rangle=2s-n.
	\end{equation}	We denote $\widehat{F}_{X_n(K)\setminus B^\circ ,\mathfrak{v}_s}:\widehat{\text{HF}}(S^3)\to \widehat{\text{HF}}(S^3_n(K)) $ and $\widehat{\nu}(K)=|\varepsilon(K)|\cdot (2\nu(K)-1)$, where $B\subseteq X_n(K)$ denotes the $0$-handle. One could refer \cite[Section 2]{HW26} for a review of the definitions of $\varepsilon(K)$ and $\nu(K)$.
	
	\begin{theorem}[{\cite[Theorem 1.10]{HW26}}]\label{HW26cal}
		There exists $s\in\mathbb{Z}$  such that $\widehat{F}_{X_n(K)\setminus B^\circ ,\mathfrak{v}_s}\ne 0$ if and only if $n\le \widehat{\nu}(K)+0^{\widehat{\nu}(K)}-1$, where we denote $0^k$ as $1$ when $k=0$ and zero for other $k$.
	\end{theorem}

	Now apply Theorem \ref{thmgeneralW}, we can easily deduce that $$\fl\left(X_n(K),U\right)\ne 0$$for any knot $K$ and any integer  $n\le \widehat{\nu}(K)+0^{\widehat{\nu}(K)}-1$, which finishes the proof of Proposition \ref{thm1.0}.

	\section{Cocore annulus cobordism map and $\fl(X_n(K),K_n,\alpha_0)$}\label{sec4}
	
	In this section, we   generalize Theorem \ref{HW26cal} to the dual knot, and then prove Theorem \ref{thm1.2} , i.e., for any nonzero integer $n$, we have \begin{equation}\label{41}
		\dim \fl(\xkn,\alpha_0)\ge \max\{2\tau(K)-n,0\},
	\end{equation}where $\alpha_0=[D^{cc}]$, recall $D^{cc}$ is the cocore disk of the $2$-handle of $X_n(K)$ with boundary $K_n$.

	Note that take the  $4$-ball $B$ in $X_n(K)$ as the $0$-handle, then consider the right-handed meridian $\mu$ in $S^3$, it stays in the complement  $S^3\setminus\nu(K)$ and the image of $-\mu$ induces dual knot $K_n$ in $S^3_n(K)$. Hence there is a product cobordism $-\mu\times I$ in $X_n(K)\setminus B^\circ$, from $-\mu$ to $K_n$.
	
	Take the canonical divides $\mathcal{A}$ on the annulus, we form a knot cobordism $$\left(X_n(K)\setminus B^\circ,  (-\mu\times I,\mathcal{A})\right):(S^3,-\mu)\to (S^3_n(K),K_n),$$where we omit the two basepoints for simplicity. We call this decorated cobordism the \textit{cocore annulus cobordism}, since one can easily see that if we fill the annulus by the disk bounded by the unknot $\mu\subset S^3$,  the new disk is isotopic to the cocore disk $D^{cc}$.  Hence the class of the Floer lasagna filling \begin{equation}
		\mathcal{F}:=\left(\{B,U,1\},(-\mu\times I,\mathcal{A})\right)
	\end{equation}lies in $\mathcal{FL}(\xkn,\alpha_0)$.
	
	We now outline the proof of \eqref{41}, we  first use the dual knot mapping cone formula developed in \cite{HL24} to identify the  knot Floer cobordism map  induced by cocore annulus cobordism with a canonical inclusion map of the mapping cone, follows the idea of \cite{OS08}. Then we can easily calculate the rank of the inclusion from simple homological algebra and hence the rank of the cocore annulus cobordism map. Now using the observation of \eqref{evin}, we will have the desired result.
	
	\subsection{Mapping cone formula of dual knot}\label{mappingconeformula}

	In \cite{HL24} , Hedden and Levine defined a mapping cone formula that for $n\ne 0$ computes $\cfk^\infty(Y_n(K),K_n)$. By restricting their mapping cone formula to $\widehat{\cfk}$, we can compute $\widehat{\cfk}(Y_n(K),K_n)$ and their hat version formula agrees with Eftekhary’s results in \cite{Eft18}.

	Let $K\subset Y$ be a rationally null-homologous knot which represents a class of order $d>0$ in $\text{H}_1(Y;\mathbb{Z})$, and $\mu$ is its right-handed meridian such that $\m{lk}(K,\mu)=1$. A relative $\spc$ structure is a homology class of nowhere-vanishing vector fields on $Y\setminus \nu(K)$ which is tangent to the boundary along $\partial\nu(K)$. We denote the set of  relative $\spc$ structure as $\rspc(Y,K)$. Given an orientation, Ozsv\'ath and Szab\'o define a map $$G_{Y,K}:\rspc(Y,K)\to \spc(Y),$$and the fibers of $G_{Y,K}$ are precisely the orbits of $\rspc(Y,K)$ under the action of $\langle\text{PD}[\mu]\rangle\subseteq \text{H}^2(Y,K;\mathbb{Z})$.

	The Alexander grading of each $\xi\in \rspc(Y,K)$ is defined as \begin{equation}\label{alex}
		A_{Y,K}(\xi):=\dfrac{\langle c_1(\xi),[F]\rangle + [\mu]\cdot [F]}{2[\mu]\cdot [F]}\in\dfrac{1}{2d}\mathbb{Z},
	\end{equation}where $[F]$ is a rational Seifert surface for $K$. Indeed, any $\xi\in\rspc(Y,K)$ is uniquely determined by the pair $(G_{Y,K}(\xi),A_{Y,K}(\xi))$ .

	In our case, we choose $Y=S^3$, hence $K$ is null-homologous, i.e., $d=1$, with Seifert framing $\lambda$.   Hence an   $n$-framing on $K$ is specified by a choice of longitude $\eta_n=\lambda+n\cdot \mu$. Let $Y_n=S^3_n(K)$ be the $n$-framed surgery on $K$, the meridian $\mu$ then is isotopic to a core circle of the glued-in solid torus. Let $K_n$ denote this core circle, with the orientation inherited from the left-handed meridian $-\mu$. 
	
	The sets $\rspc(Y,K)$ and $\rspc(Y_n,K_n)$  are canonically identified. And the orientation of $K_n$ also induces a map $G_{Y_n,K_n}:\rspc(Y_n,K_n)\to \spc(Y_n)$ whose fibers are the orbits of the action of $\langle\text{PD}[\eta_n]\rangle$.

	Assume henceforth that $n\ne 0$ in this section. Choose a $\spc$ structure $\mathfrak{t}$ on $Y_n(K)$. Index the elements of $G_{Y_n,K_n}^{-1}(\mathfrak{t})$ by $(\xi_l)_{l\in\mathbb{Z}}$, where $\xi_{l+1}=\xi_l+\text{PD}[\eta_n]$. Since $Y=S^3$ has unique $\spc$ structure $\mathfrak{s}_0$, then we know that for any $l\in\mathbb{Z}$, we have $G_{S^3,K}(\xi_l)=\mathfrak{s}_0$. Denote $s_l=A_{S^3,K}(\xi_l)\in\mathbb{Z}$, we obtain\begin{equation}
		s_{l+1}-s_l=A_{S^3,K}(\xi_{l+1})-A_{S^3,K}(\xi_l)=n
	\end{equation}from \eqref{alex}. Moreover, it is easy to see that \begin{equation}\label{xilalex}
		A_{Y_n,K_n}(\xi_l)=\dfrac{2A_{S^3,K}(\xi_l)+n-1}{2n}=\dfrac{2s_l+n-1}{2n}.
	\end{equation}
	
	\begin{remark}
Note that any $\xi \in \rspc(S^3, K) \allowbreak \cong \allowbreak \rspc(S^3_n(K), K_n)$ is determined by $G_{S^3,K}(\xi) = \mathfrak{s}_0$ and $A_{S^3,K}(\xi) = s \in \mathbb{Z}$. Hence there is a one-to-one correspondence among $\rspc(S^3_n(K),K_n)$, $\rspc(S^3, K)$, and $\mathbb{Z}$.
	\end{remark}

	Now we can state the dual knot mapping cone formula. For simplicity, we only state their hat version here. Let $C=\cfk^\infty(S^3,K)$ with $i,j$ filtration be the knot Floer complex defined in \cite{OS04c}, which is generated by the elements of $[\textbf{x},i,j]$ where $\textbf{x}\in \mathbb{T}_{\bm{\alpha}}\cap \mathbb{T}_{\bm{\beta}}$. And $C$ admits a $\mathcal{U}$ action defined by $\mathcal{U}[\textbf{x},i,j]=[\textbf{x},i-1,j-1]$.
	
	\begin{theorem}[{\cite[Proposition 1.5]{Eft18}}, {\cite[Corollary 3.6]{HL24}}]
		Let $\xi \in \rspc(S^3_n(K),K_n) \allowbreak \cong \allowbreak \rspc(S^3, K)$ and set $s = A_{S^3,K}(\xi)$. Then $\widehat{\hfk}(S^3_n(K),K_n,\xi)$ is isomorphic to  the homology of the mapping cone of \begin{equation}\label{hvs}
			(h_{s},v_{s+n}): C\{i\le 0, j=s\}\oplus C\{i=0,j\le s+n-1 \} \to C\{i=0\}, 
		\end{equation}where $h_{s}$ is the composition of multiplication of  $\mathcal{U}^{s}$ with  the standard flip map interchanging the $i,j$ filtrations, $v_{s+n}$ is the inclusion. And we denote this mapping cone as $\mathbb{X}_{s}$.
	\end{theorem}

	For each $s\in \mathbb{Z}$, there is a  $\mathfrak{u}_s\in\spc(X_n(K))$ such that \begin{equation}\label{uspin}
		\left\langle c_1(\mathfrak{u}_s),[\widehat{F}]  \right\rangle=-2s-n.
	\end{equation}
	
	\begin{lemma}\label{xishift}
		The knot Floer cobordism map induced from the cocore annulus cobordism satisfies $$\widehat{F}_{X_n(K)\setminus B^\circ,(-\mu\times I,\mathcal{A}),\mathfrak{u}_{s}}:\widehat{\hfk}(S^3,-\mu)\to \widehat{\hfk}(S^3_n(K),K_n,\xi).$$
	\end{lemma}
	
	\begin{proof}
		It suffices to calculate the Alexander grading shift of the  map $\widehat{F}_{X_n(K)\setminus B^\circ,(-\mu\times I,\mathcal{A}),\mathfrak{u}_{s}}$.
		Its calculation is similar to that in \cite[Lemma 3.2]{HHSZ26}, with minor notational differences.
		Let $W=X_n(K)\setminus B^\circ$, $\Sigma_\mu=-\mu\times I$. By choosing a smaller  $B$ , we can assume $\widehat{F}$ lies in the interior of $W$, hence it is easy to see $\m{H}_2(W;\mathbb{Q})\cong \m{H}_2(W,\partial W;\mathbb{Q})$ are all generated by $[\widehat{F}]$. And $[\widehat{F}]\cdot [\widehat{F}]=n$.
		
		We view $[\Sigma_\mu]\in \m{H}_2(W,\partial W;\mathbb{Q})$ by capping the disk of $-\mu\times \{0\}$  in $S^3=\partial B$, and take $[\widehat{\Sigma}_\mu]$ as the preimage of $[\Sigma_\mu]$ under the natural isomorphism between $\m{H}_2(W;\mathbb{Q})$ and $\m{H}_2(W,\partial W;\mathbb{Q})$. Easy to see $[\Sigma_\mu]\cdot [\widehat{F}]=[\widehat{\Sigma}_\mu]\cdot [\widehat{F}]=-1$, this is because we can view the intersection as occurring between $K$ and a disk $D$ in $S^3$ such that $\partial D=-\mu$, now we have $[\widehat{\Sigma}_\mu]\cdot [\widehat{F}]=\m{lk}(K,-\mu)=-1$. Hence from $[\widehat{F}]\cdot [\widehat{F}]=n$ we have $$[\Sigma_\mu]=-\dfrac{1}{n}[\widehat{F}]\in \m{H}_2(W,\partial W;\mathbb{Q}),\quad [\widehat{\Sigma}_\mu]=-\dfrac{1}{n}[\widehat{F}]\in \m{H}_2(W;\mathbb{Q}).$$

		Now from Proposition  \ref{alexandershift} and \eqref{uspin}, we have the Alexander grading shift of $\widehat{F}_{X_n(K)\setminus B^\circ,(-\mu\times I,\mathcal{A}),\mathfrak{u}_{s}}$ is \begin{equation*}
			\dfrac{\langle c_1(\mathfrak{u}_s), [\widehat{\Sigma}_\mu]\rangle -[\Sigma_\mu]\cdot [\widehat{\Sigma}_\mu]}{2}=\dfrac{1}{2}\left(-\dfrac{1}{n}\cdot (-2s-n)-\dfrac{1}{n^2}\cdot n\right)=\dfrac{2s+n-1}{2n}.
		\end{equation*}Hence from $\widehat{\hfk}(S^3,\mu)$ supports only in Alexander grading zero, we have  the target of $\widehat{F}_{X_n(K)\setminus B^\circ,(-\mu\times I,\mathcal{A}),\mathfrak{u}_{s}}$ has Alexander grading $\dfrac{2s+n-1}{2n}=A_{S^3_n(K),K_n}(\xi)$ from \eqref{xilalex}.  Now since $\xi$ is uniquely determined by $s=A_{S^3,K}(\xi)$,  we know the result holds.
	\end{proof}

	Now we can state our main result in this section,

	\begin{theorem}\label{inclusion}
		Consider the canonical inclusion  $\iota_{s}:C\{i=0\}\to \mathbb{X}_{s}$, then its induced  map on  homology $$\iota_{s,*}:\mathrm{H}(C\{i=0\})\to \mathrm{H}(\mathbb{X}_s)$$ can be identified with the cocore annulus cobordism map $\widehat{F}_{X_n(K)\setminus B^\circ,(-\mu\times I,\mathcal{A}),\mathfrak{u}_{s}}$.
	\end{theorem}

	Note that $\text{H}(C\{i=0\})\cong \mathbb{F}\cong \widehat{\hfk}(S^3,-\mu)$, hence from the theorem above, we know that $$\text{rank}\ \widehat{F}_{X_n(K)\setminus B^\circ,(-\mu\times I,\mathcal{A}),\mathfrak{u}_{s}}=\text{rank}
	\ \iota_{s,*}\in\{0,1\}.  $$ In particular, from the easy algebraic lemma below, we can determine for which $s$, the rank is nonzero.

	\begin{lemma}
		Consider the chain map $f=(f_1,f_2):A_1\oplus A_2\to B$, and the mapping cone $\mathbb{X}=\mathrm{Cone}(f)$, together with the inclusion $\iota: B\to \mathbb{X}$. If $\text{H}(B)\cong \mathbb{F}$, then we have $\iota_*\ne 0$ if and only if  $f_*=0$, i.e., $f_{1,*}=0$ and $f_{2,*}=0$. 
	\end{lemma}
	
	\begin{proof}
	Consider the associated exact triangle of the mapping cone $\mathbb{X}=\mathrm{Cone}(f)$, that is, $$\text{H}(A_1)\oplus \text{H}(A_2)\xrightarrow{f_*=(f_{1,*},f_{2,*})}{} \text{H}(B)\overset{\iota_*}{\longrightarrow}\text{H}(\mathbb{X})\to \text{H}(A_1)\oplus \text{H}(A_2).$$
		
		Since $\text{H}(B)\cong \mathbb{F}$, then we have $\iota_*$ is either trivial or injective. Hence $\iota_*\ne 0$ is equivalent to $\iota_*$ is injective, i.e., $\ker \iota_*=\{0\}$, then $\m{Im} f_*=\ker \iota_*=\{0\} $, i.e., $f_*=0$.
	\end{proof}

	Hence $\text{rank}\  \iota_{s,*}=1$ if and only if $h_{s,*}=0$ and $v_{s+n,*}=0$. Recall \begin{equation}
		\tau(K)=\min\{s-1|v_{s}:C\{i=0,j\le s-1\}\to C\{i=0\}\text{ induces surjection on homology} \},
	\end{equation}hence $v_{s,*}= 0$ if and only if $s-1<\tau(K)$, i.e., $s\le \tau(K)$.

	Note that $h_s$ is the composition of the following maps $$C\{i\le 0,j=s\}\xrightarrow{\mathcal{U}^{s}} C\{i\le -s,j=0\}\xrightarrow{\simeq} C\{i=0,j\le -s\}\xrightarrow{v_{-s+1}}C\{i=0\},$$where the first two maps are homotopy equivalences, hence $h_{s,*}= 0$ if and only if $v_{-s+1,*}= 0$, i.e., $-s+1-1<\tau(K)$, hence $s\ge 1-\tau(K)$.

	In summary, if $v_{s+n,*}=0$ and $h_{s,*}=0$, then we have \begin{equation}
		1-\tau(K)\le s\le \tau(K)-n,
	\end{equation}and for any other $s$, either $v_{s+n,*}\ne 0$ or $h_{s,*}\ne 0$.

	Now combine the result above, we have 
	
	\begin{corollary}\label{rankcocore}
		For any integer $s$, $\m{rank}\ \widehat{F}_{X_n(K)\setminus B^\circ,(-\mu\times I,\mathcal{A}),\mathfrak{u}_{s}}=1$ if and only if $1-\tau(K)\le s\le \tau(K)-n$, otherwise the rank is zero. Hence consider the map $$\bigoplus_{s\in\mathbb{Z}} \widehat{F}_{X_n(K)\setminus B^\circ,(-\mu\times I,\mathcal{A}),\mathfrak{u}_{s}}:\bigoplus_{s\in\mathbb{Z}}\widehat{\hfk}(S^3,-\mu)\to \widehat{\hfk}(S^3_n(K),K_n),$$we have \begin{equation}
			\m{rank}\left(\bigoplus_{s\in\mathbb{Z}}\widehat{F}_{X_n(K)\setminus B^\circ,(-\mu\times I,\mathcal{A}),\mathfrak{u}_{s}}\right)=\max\{2\tau(K)-n,0\}.
		\end{equation}
	\end{corollary}
	
	\subsection{Proof of Theorem \ref{thm1.2}} Now we  first use Corollary \ref{rankcocore} to give a proof of Theorem \ref{thm1.2}, and  defer the proof of Theorem \ref{inclusion} to Section  \ref{holomo}.

	Let \(K_+\) and \(K_-\) be oppositely oriented \(n\)-framed parallels of \(K\) and the orientation of $K_+$ is the same as $K$. Now consider $L_1=-\mu\cup K_+\cup K_-$ and each has two basepoints. Let $\Sigma= (-\mu\times I)\cup D_+\cup D_-$, where $D_\pm$ are the core disks bounded by $K_\pm$, it is easy to see the capped off disks $\widehat{D}_{\pm}$ satisfying $[\widehat{D}_\pm]=\pm [\widehat{F}]$, hence we have the capped off surface $[\widehat{\Sigma}]=[D^{cc}]+[\widehat{F}]-[\widehat{F}]=[D^{cc}]=\alpha_0$ . Choose canonical divides on $-\mu\times I$ as before and also canonical divides  on disks $D_\pm$, we have a decorated surface $(\Sigma,\mathcal{A}_1)$.

	\begin{remark}
		Let $K_\pm$ live in a smaller tubular neighborhood of $K$, and $\mu$ lies in the boundary of a larger tubular neighborhood, it is easy  to ensure two core disks $D_\pm$ are disjoint from the annulus $-\mu\times I$. 
	\end{remark}
	
	Now furthermore we push the incoming \(S^3\) slightly into the \(0\)-handle,  and we can consider the pair-of-pants cobordism $\mathcal{P}$ in $S^3\times [1-\varepsilon,1]$, from $-\mu\cup U$ to $L_1=-\mu\cup K_+\cup K_-$, which is defined by Chen in \cite{chen2022floer}. More precisely $\mathcal{P}$ is the product annulus of $-\mu$ together with a pair-of-pants from $U$ to $K_+\cup K_-$.

	\begin{figure}[htbp]
		\centering
		\includegraphics[width=0.6\linewidth]{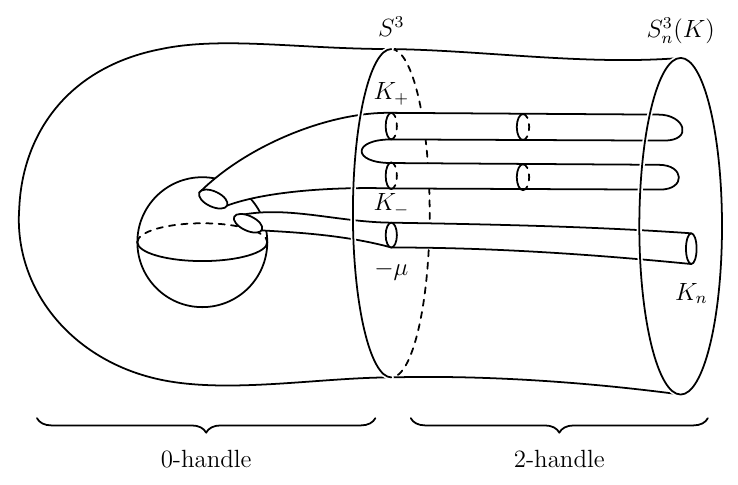}
		\caption{A schematic picture of cobordism $(\Sigma,\mathcal{A}_1)\circ \mathcal{P}$ where the dividing curves are omitted.}
		\label{fig:handleannulustikz}
	\end{figure}

	Note that $$\widehat{\hfl}(-\mu\cup U)\cong \widehat{\hfl}(-\mu)\otimes V,$$where $V=\mathbb{F}\langle T,B\rangle$ where for the death map $\mathcal{D}$ of $U$, \begin{equation}
		\mathcal{D}(T)=0,\quad \mathcal{D}(B)=1.
	\end{equation}

	Now consider $x=\widehat{F}_\mathcal{P}(1\otimes B)\in\widehat{\hfl}(L_1)= \widehat{\hfl}(-\mu\cup K_+\cup K_-)$. Note that $(\Sigma,\mathcal{A}_1)\circ \mathcal{P}$ is isotopic to $(-\mu\times I)\cup \mathcal{D}$ in $X_n(K)\setminus B^\circ$ from $(S^3,-\mu\cup U)$ to $(S^3_n(K),K_n)$. Hence from the composition law \eqref{linkfloercompo}, we have \begin{align*}
		\widehat{F}_{X_n(K)\setminus B^\circ, (\Sigma,\mathcal{A}_1),\mathfrak{u}_s}(x)&=\widehat{F}_{X_n(K)\setminus B^\circ, (\Sigma,\mathcal{A}_1),\mathfrak{u}_s}\circ \widehat{F}_\mathcal{P}(1\otimes B)\\
		&=\widehat{F}_{X_n(K)\setminus B^\circ, (-\mu\times I,\mathcal{A})\cup \mathcal{D},\mathfrak{u}_s}(1\otimes B)\\
		&=\widehat{F}_{X_n(K)\setminus B^\circ, (-\mu\times I,\mathcal{A}),\mathfrak{u}_s}(1)\otimes \mathcal{D}(B)\\
		&=\widehat{F}_{X_n(K)\setminus B^\circ, (-\mu\times I,\mathcal{A}),\mathfrak{u}_s}(1)\in \widehat{\hfk}(S^3_n(K),K_n,\xi),
	\end{align*}where $\xi$ is determined by $A_{S^3,K}(\xi)=s$.

	Now note that on $\widehat{\hfl}(L_1)= \widehat{\hfl}(-\mu\cup K_+\cup K_-)$, we have an Alexander multi-grading  $(A_\mu,A_+,A_-)$. Now for a homogeneous element $u$, if $\widehat{F}_{X_n(K)\setminus B^\circ, (\Sigma,\mathcal{A}_1),\mathfrak{u}_s}(u)\ne 0$, we can  calculate the $A_\pm$ gradings of $u$ by calculating the Alexander grading shift of the cobordism map. Note that $\Sigma=\Sigma_\mu\cup D_+\cup D_-$, where $\Sigma_\mu=-\mu\times I$ as the notation in the proof of Lemma \ref{xishift}, then $[\Sigma]=[\Sigma_\mu]=-\dfrac{1}{n}[\widehat{F}]$, $[\widehat{D}_\pm]=\pm [\widehat{F}]$. It follows from Proposition \ref{alexandershift} and  Remark \ref{emptyremark}, we have \begin{align*}
		0-A_+(u)&=\dfrac{\langle c_1(\mathfrak{u}_s),[\widehat{D}_+]\rangle-[\Sigma]\cdot  [\widehat{D}_+]}{2}=\dfrac{-2s-n+1}{2},\\
		0-A_-(u)&=\dfrac{\langle c_1(\mathfrak{u}_s),[\widehat{D}_-]\rangle-[\Sigma]\cdot  [\widehat{D}_-]}{2}=\dfrac{2s+n-1}{2}.
	\end{align*}

	Hence we can decompose $x$ by the Alexander grading $A_+$, i.e., $$x=\sum_{s\in\mathbb{Z}}x_s,\quad A_+(x_s)=\dfrac{2s+n-1}{2}.$$ Now from the Alexander grading, we can easily deduce \begin{equation*}
		\widehat{F}_{X_n(K)\setminus B^\circ, (\Sigma,\mathcal{A}_1),\mathfrak{u}_s}(x_t)=\begin{cases}
			\widehat{F}_{X_n(K)\setminus B^\circ, (\Sigma,\mathcal{A}_1),\mathfrak{u}_s}(x) & \text{ if } t=s \\
			0& \text{ if } t\ne s
		\end{cases}.
	\end{equation*}
	
	Take $\displaystyle \widehat{F}_{X_n(K)\setminus B^\circ, (\Sigma,\mathcal{A}_1)}=\sum_{s\in\mathbb{Z}}\widehat{F}_{X_n(K)\setminus B^\circ, (\Sigma,\mathcal{A}_1),\mathfrak{u}_s}$, then we have \begin{equation*}
		\widehat{F}_{X_n(K)\setminus B^\circ, (\Sigma,\mathcal{A}_1)}(x_s)=\widehat{F}_{X_n(K)\setminus B^\circ, (-\mu\times I,\mathcal{A}),\mathfrak{u}_s}(1)\in \widehat{\hfk}(S^3_n(K),K_n,\xi),
	\end{equation*}hence for $s\ne s'$, then $A_{S^3_n(K),K_n}(\xi)\ne A_{S^3_n(K),K_n}(\xi') $ from \eqref{xilalex}, we have \begin{equation}\label{ssprime}
		\widehat{F}_{X_n(K)\setminus B^\circ, (\Sigma,\mathcal{A}_1)}(x_s)\ne \widehat{F}_{X_n(K)\setminus B^\circ, (\Sigma,\mathcal{A}_1)}(x_{s'})
	\end{equation}if both of them are nonzero.
	
	Now we take Floer lasagna filling \begin{equation}
		\mathcal{F}_s:=\left(\{B,L_1,x_s\},(\Sigma,\mathcal{A}_1)\right),
	\end{equation}then from the definition of evaluation map \eqref{ev}, we have \begin{equation*}
		\m{ev}\left([\mathcal{F}_s]\right)=	\widehat{F}_{X_n(K)\setminus B^\circ, (\Sigma,\mathcal{A}_1)}(x_s)=\widehat{F}_{X_n(K)\setminus B^\circ, (-\mu\times I,\mathcal{A}),\mathfrak{u}_s}(1).
	\end{equation*}

	Now from Corollary \ref{rankcocore}, for any $1-\tau(K)\le s\le \tau(K)-n$, we have $\widehat{F}_{X_n(K)\setminus B^\circ, (-\mu\times I,\mathcal{A}),\mathfrak{u}_s}(1)\ne 0$ hence $\m{ev}\left([\mathcal{F}_s]\right)$ is nonzero, and for different $1-\tau(K)\le s,s'\le \tau(K)-n$, from \eqref{ssprime} we have $\m{ev}\left([\mathcal{F}_s]\right)\ne \m{ev}\left([\mathcal{F}_{s'}]\right)$.
	
	Thus in summary, we have $$\dim \fl(X_n(K),K_n,\alpha_0)\ge \mathrm{rank} \text{ ev}\ge \max\{2\tau(K)-n,0\}$$which is the desired result.
	
	\subsection{Holomorphic triangle and the inclusion map}\label{holomo}In this section, we  prove Theorem \ref{inclusion}.

	Let $K\subset S^3$ be an oriented knot, with doubly based Heegaard diagram $$\left(\Sigma, \bm{\alpha}=\{\alpha_i\}_{i=1}^g, \bm{\beta}=\{\beta_i\}_{i=1}^g,w,z\right)$$such that $\beta_g$ is isotopic to the meridian of $K$, and $\beta_g$ meets $\alpha_g$ in a single point  and is disjoint from the remaining $\alpha$ curves. There is an arc $t_\alpha$ from $z$ to $w$ that meets $\beta_g$ in a single point and is disjoint from all other $\alpha$ and $\beta$  curves.

Next, consider the curve $\gamma_g$, which is an \(n\)-framed longitude meeting \(\beta_g\) once and disjoint from the remaining \(\beta\) curves. It is oriented in the same direction as \(K\), and \(w, z\) lie to its right. For \(i=1,\cdots,g-1\), let \(\gamma_i\) be a small pushoff of \(\beta_i\) that intersects \(\beta_i\) in two points.

Taking a tubular neighborhood of \(\beta_g\), we can choose another basepoint \(z'\) as shown in Figure \ref{fig:heegaard-dia} below.

	\begin{figure}[htbp]
		\centering
		\includegraphics[width=0.6\linewidth]{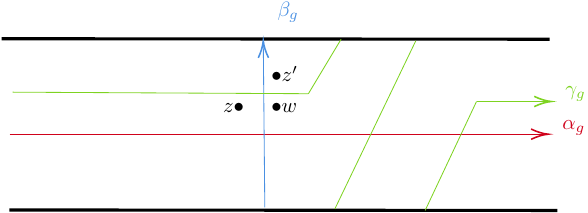}
		\caption{tubular neighborhood of $\beta_g$ in $\Sigma$}
		\label{fig:heegaard-dia}
	\end{figure}

	Hence we have a Heegaard triple $$\mathcal{T}=\left(\Sigma,\bm{\alpha}=\{\alpha_i\}_{i=1}^g, \bm{\beta}=\{\beta_i\}_{i=1}^g, \bm{\gamma}=\{\gamma_i\}_{i=1}^g ,w,z'\right),$$such that \begin{align*}
		\left(\Sigma,\bm{\alpha}, \bm{\beta},w,z'\right)&=(S^3,-\mu)\\
		\left(\Sigma, \bm{\beta},\bm{\gamma},w,z'\right)&=\left(\#^{g-1}(S^1\times S^2),U\right)\\
		\left(\Sigma,\bm{\alpha}, \bm{\gamma},w,z'\right)&=(S^3_n(K),K_n).
	\end{align*}

	From \cite[Section 3.1]{Zem19b}, we have a link cobordism $(W_{\bm{\alpha\beta\gamma}}, \Sigma_{\bm{\alpha\beta\gamma}})$ from a doubly based Heegaard triple $\mathcal{T}$. Easy to see the underlying  $4$-manifold is $X_n(K)$ with a $4$-ball $B$ and $\natural^{g-1}(S^1\times \mathbb{D}^3)$ removed (as in \cite[Example 8.1]{OS04a}) and the link cobordism is a pair of pants . Hence after filling $Y_{\bm{\beta\gamma}}$ , we get $X_n(K)\setminus B^\circ$, we can also use a disk to cap off the unknot $U$ in $Y_{\bm{\beta\gamma}}=\#^{g-1}(S^1\times S^2)$,and this cobordism   sends $1$ to the top-generator $\Theta_{\bm{\beta\gamma}}^+$.
	
	Now \cite[Theorem 9.3]{JZ24} shows that the link Floer cobordism map induced from the Heegaard triple is the same as the holomorphic  triangle counting map. Composing the triple cobordism with  the birth map, we obtain $$F_{X_n(K)\setminus B^\circ,(-\mu\times I,\mathcal{A}),\mathfrak{u}_{s}}:\widehat{\cfk}(S^3,-\mu)\to \widehat{\cfk}(S^3_n(K),K_n,\xi)$$is chain homotopic to the following holomorphic triangle counting map\begin{equation}
		F_{\mathcal{T},\mathfrak{u}_s}(\textbf{x})=\sum_{\textbf{y}\in\mathbb{T}_\alpha\cap \mathbb{T}_\gamma}\sum_{\substack{\psi \in \pi_2(\textbf{x},\Theta_{\bm{\beta\gamma}}^+,\textbf{y}) \\ \mu(\psi) = 0, n_w(\psi)=n_{z'}(\psi)=0\\ \mathfrak{s}_w(\psi)=\mathfrak{u}_s}} 
		\# \mathcal{M}(\psi) \cdot \textbf{y}.
	\end{equation}

	Now recall in \cite{Eft18}, for any $\xi\in \rspc(S^3,K)=\rspc(S^3_n(K),K_n)$,  Eftekhary constructed a triangle 
	
	\[\begin{tikzcd}
		{\widehat{\text{CF}}(\bm{\alpha},\bm{\beta},w,z')} && {\widehat{\text{CF}}(\Sigma,\bm{\alpha},\bm{\gamma},w,z',\xi)} \\
		& {C\{i\le 0, j=s\}\oplus C\{i=0,j\le s+n-1 \} }
		\arrow["{h^\xi}", from=1-1, to=1-3]
		\arrow["{g^\xi}", from=1-3, to=2-2]
		\arrow["{f^\xi}", from=2-2, to=1-1]
	\end{tikzcd}\]where the definition of $h^\xi$ is exactly 
	same as our $F_{\mathcal{T},\mathfrak{u}_s}$, and 	we need to modify our $f^\xi,g^\xi$ from \cite{Eft18} by composing the identification between $C\{i\le 0, j=s\}\oplus C\{i=0,j\le s+n-1 \}$ and Eftekhary's $C_{n,m}(\xi)$. Under this modification,  $f^\xi$ agrees with our $(h_s,v_{s+n})$ defined in \eqref{hvs}.

	Moreover, there is a  chain homotopy \(H^\xi\) from \(f^\xi\circ g^\xi\) to zero, then \cite[Theorem 2.3]{Eft18} proves that \begin{equation}
		\jmath^\xi=(g^\xi,H^\xi): \widehat{\text{CF}}(\Sigma,\bm{\alpha},\bm{\gamma},w,z',\xi)\to \text{Cone}(f^\xi)=\mathbb{X}_s
	\end{equation}is a quasi-isomorphism. 
	
	For short, we denote 
	\begin{itemize}
		\item $B=\widehat{\text{CF}}(\bm{\alpha},\bm{\beta},w,z')$,
		\item $A_s=C\{i\le 0, j=s\}\oplus C\{i=0,j\le s+n-1 \}$,
	\end{itemize}it is easy to see $\m{H}(B)\cong \widehat{\hfk}(S^3,U)\cong \mathbb{F}$. And we have  the following exact triangle of the mapping cone, \[\begin{tikzcd}
		{\m{H}(B)} && {\m{H}(\mathbb{X}_s)} \\
		& {\m{H}(A_s) }
		\arrow["{\iota_{s,*}}", from=1-1, to=1-3]
		\arrow["{p_*}", from=1-3, to=2-2]
		\arrow["{f^\xi_*}", from=2-2, to=1-1]
	\end{tikzcd}\]

	Note that   $p\circ \jmath^\xi=g^\xi$ from the definition of $\jmath^\xi$. Since $\jmath^\xi$ is a quasi-isomorphism, we have $$\m{Im} \left(\jmath^\xi_*\circ h^\xi_*\right)=\jmath_*^\xi\m{Im } h^\xi_*= \jmath_*^\xi\ker g^\xi_*=\ker p_*=\m{Im } \iota_{s,*},$$then from $\m{H}(B)\cong \mathbb{F}$ is one dimensional, we have $\jmath^\xi_*\circ h^\xi_*=\iota_{s,*}$. Hence we have the following commutative diagram, \[\begin{tikzcd}
		{\widehat{\text{HFK}}(S^3,-\mu)} & {} & {\widehat{\text{HFK}}(S^3_n(K),K_n,\xi)} \\
		{\text{H}(C\{i=0\})} & {} & {\text{H}(\mathbb{X}_s)} \\
		{}
		\arrow["{h^\xi_*=F_{\mathcal{T},\mathfrak{u}_s,*}}", from=1-1, to=1-3]
		\arrow["{=}"', from=1-1, to=2-1]
		\arrow["{\jmath^\xi_*}", from=1-3, to=2-3]
		\arrow["{\iota_{s,*}}"', from=2-1, to=2-3]
	\end{tikzcd}\]Note that $\jmath^\xi_*$ is an isomorphism,  we finish the proof of Theorem \ref{inclusion} and hence the proof of Theorem \ref{thm1.2}.

	\section{Surgery exact triangle and $\fl(X_n(K),K_n,\alpha_{-k})$ for $k\ge 0$ }\label{sec5}

	In this section, we  use the knot Floer surgery exact triangle to prove some nonvanishing results of $\fl(X_n(K),K_n,\alpha_{-k})$ for $k\ge 0$.
	
	\subsection{Surgery exact triangle in knot Floer homology}We first state a knot Floer version of the surgery exact triangle.
	
	\begin{theorem}\label{knotexact}
		Let \(K\subset S^3\) be a knot, and let \(J\subset S^3\setminus\nu(K)\) be another knot. For any integer $n$, we have a surgery exact triangle \begin{equation}\label{knotfloertriangle}
			\begin{tikzcd}
				{\widehat{\hfk}(S^3,J)} \arrow[rr, "F_1"] &                                                         & {\widehat{\hfk}(S^3_n(K),J)} \arrow[ld, "F_2"] \\
				& {\widehat{\hfk}(S^3_{n+1}(K),J )}\arrow[lu, "F_3"]                         
			\end{tikzcd}
		\end{equation}where $F_i$ is the sum, over all $\spc$ structures, of the knot Floer cobordism maps induced by the surgery cobordisms with a product annulus traced by $J$. 
	\end{theorem}
	
	\begin{proof}
		Suppose \((M, \gamma)\) is a balanced sutured manifold and \(K\) is a knot in \(M\). Consider three balanced sutured manifolds \((M_i, \gamma_i)\) for \(i = 1, 2, 3\) obtained from \((M, \gamma)\) by Dehn surgeries along \(K\). If the Dehn filling curves \(\eta_1, \eta_2, \eta_3 \subset \partial\nu(K)\) satisfy
		\[
		\eta_1 \cdot \eta_2 = \eta_2 \cdot \eta_3 = \eta_3 \cdot \eta_1 = -1,
		\]then in \cite[Theorem 4.14]{LY24}, they showed that there is an exact triangle 	\begin{equation}\label{sutureexacttriangle}
			\begin{tikzcd}
				{SFH(M_1,\gamma_1)} \arrow[rr, "F_1"] &                                       & {SFH(M_2,\gamma_2)} \arrow[ld, "F_2"] \\
				& {SFH(M_3,\gamma_3)} \arrow[lu, "F_3"] &                                      
			\end{tikzcd}
		\end{equation}where the map $F_i$ is the sutured Floer cobordism map.

		Note that  {\cite[Proposition 9.2]{Juh06}} showed that $SFH(M,\gamma)\cong \widehat{\hfk}(Y,J)$ for any $M=Y\setminus \nu(J)$ and the suture $\gamma$ is two parallel meridians of $J$ with opposite orientations. This identifies sutured Floer homology with knot Floer homology. In \cite{JZ21}, they further identified the above sutured Floer cobordism map and hat version link
		Floer cobordism map when the sutured manifold is link complement. Hence under these identifications, one can easily deduce the desired exact triangle \eqref{knotfloertriangle} from the sutured Floer exact triangle \eqref{sutureexacttriangle}.\end{proof}

	\begin{remark}
		Cai Zeng pointed out to the author that an alternative proof of the surgery exact triangle above can be obtained as follows: apply the bordered sutured Floer functor
		$\widehat{\mathit{\m{BSA}}}(S^3\setminus \nu(K),J)\boxtimes \cdot$
		to the exact triangle of bordered Floer type \(D\) structures relating
		\(\widehat{\mathit{\m{CFD}}}(\mathcal{H}_\infty)\),
		\(\widehat{\mathit{\m{CFD}}}(\mathcal{H}_n)\), and
		\(\widehat{\mathit{\m{CFD}}}(\mathcal{H}_{n+1})\),
		then invoke the bordered-sutured pairing theorem and take homology.
	\end{remark}
	
	Now to prove the nonvanishing result of $\fl(X_n(K),K_n,\alpha_{-k})$, we want to find a Floer lasagna filling representative, such that its image under the evaluation map is nonzero. We realize the evaluation of a suitable filling as one of the maps in a surgery exact triangle.
	
	\subsection{A special surgery triad }\label{surgerytraid} Let $K\subset S^3$ be an oriented knot with tubular neighborhood $\nu(K)$, and $\mu, \lambda\in \partial\nu(K)$ be its right-handed meridian and Seifert longitude. Hence $\lambda\cdot \mu=1$, $\mu\cdot\lambda=-1$. We choose the Seifert framing as the zero framing. We now fix the notation as the $(p,q)$-cable is the curve $p\lambda+q\mu$ on $\partial\nu(K)$ and  $\dfrac{p}{q}$ surgery is Dehn filling along the curve $p\mu+q\lambda$. And we use $-K$ to denote $K$ with opposite orientation and $m(K)$ to denote the mirror of $K$.
	
	Let $J_{n}^{(k)}\subset \partial\nu(K)$ be a curve representing $-(kn+1)\mu-k\lambda$. Then in $S^3$, $J_n^{(k)}$ is isotopic to the cable knot $-K_{k,kn+1}$. Now consider the $n$-surgery on $K$, note that $\nu(K_n)\cong \nu (K)$ with right-handed meridian $\mu_n$ of $K_n$ is isotopic to $n\mu+\lambda$, the longitude $\lambda_n$ is isotopic to $-\mu$, hence $$[J_n^{(k)}]=[-(kn+1)\mu-k\lambda]=[-k\mu_n+\lambda_n]\in \text{H}_1(\partial\nu(K_n)).$$Hence in $S_n^3(K)$, $J_n^{(k)}$ is isotopic to $\lambda_n$, i.e., the dual knot $K_n$.
	
	Similarly, in $S^3_{n+1}(K)$, note that $\mu_{n+1}=(n+1)\mu+\lambda$, and $\lambda_{n+1}=-\mu$, hence $$[J_n^{(k)}]=[-(kn+1)\mu-k\lambda]=[-k\mu_{n+1}-(k-1)\lambda_{n+1}]\in \text{H}_1(\partial\nu(K_{n+1})),$$then we know that $J_n^{(k)}$ is isotopic to the cable knot $-(K_{n+1})_{k-1,k}$ of $K_{n+1}$.
	
	\begin{figure}[htbp]
		\centering
		\includegraphics[width=0.6\linewidth]{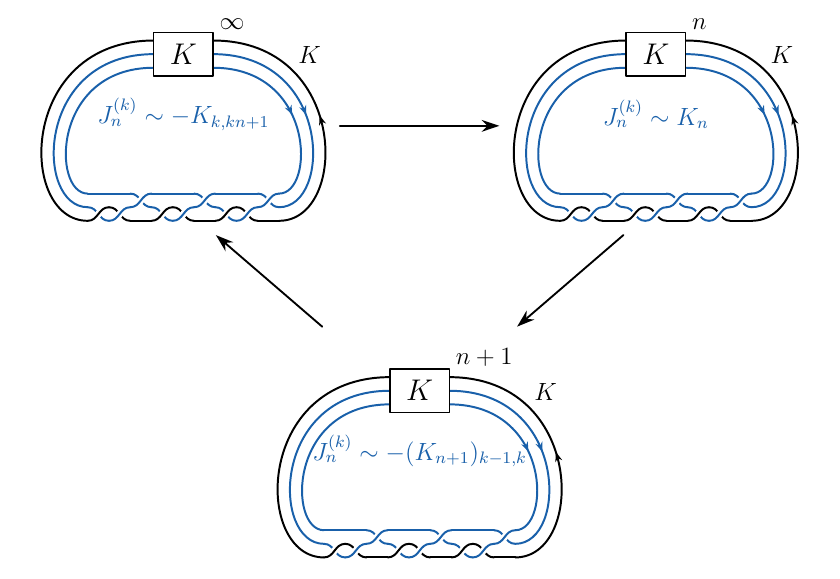}
		\caption{surgery triad schematic diagram  for $k=2$ and $n=1$ }
		\label{fig:surgerytriad-1}
	\end{figure}

	In the cobordism between \(S^3\) and \(S^3_n(K)\), tracking the image of \(J_n^{(k)}\) yields an  annulus between \(-K_{k,kn+1}\) and the dual knot \(K_n\). We denote this annulus by \(A_n^{(k)}\) and equip it with the canonical divides \(\mathcal{A}\).  Then we form a knot cobordism $$\mathcal{A}_n^{(k)}=(X_n(K)\setminus B^\circ, (A_n^{(k)},\mathcal{A})):(S^3,-K_{k,kn+1})\to (S^3_n(K),K_n).$$

	Now apply Theorem \ref{knotexact} for surgery knot $K$ and $J=J_n^{(k)}\subset S^3\setminus \nu(K)$, we have an exact triangle \begin{equation}\label{trianglecable}
		\begin{tikzcd}
			{\widehat{\hfk}(S^3,-K_{k,kn+1})} \arrow[rr, "\widehat{F}_{\mathcal{A}_n^{(k)}}"] &                                                         & {\widehat{\hfk}(S^3_n(K),K_n)} \arrow[ld, "\widehat{F}_2"] \\
			& {\widehat{\hfk}(S^3_{n+1}(K),-(K_{n+1})_{k-1,k} )}\arrow[lu, "\widehat{F}_3"] &                                                
		\end{tikzcd}
	\end{equation}where $\widehat{F}_{\mathcal{A}_n^{(k)}}$ is the knot Floer cobordism map induced from $\mathcal{A}_n^{(k)}$. We can calculate the rank of $\widehat{F}_{\mathcal{A}_n^{(k)}}$ from the dimensions of three knot Floer homology groups using the simple lemma below.

	\begin{lemma}\label{alglem}
		Suppose we have an exact triangle between three finite-dimensional $\mathbb{F}$-vector spaces, \[	\begin{tikzcd}
			{A} \arrow[rr, "f"] &                                                         & {B} \arrow[ld, "g"] \\
			& {C} \arrow[lu, "h"] &                                                
		\end{tikzcd}\]
		then we have \begin{equation}
			\mathrm{rank} f=\dfrac{\dim A+\dim B-\dim C}{2}.
		\end{equation}
	\end{lemma}
	
	\begin{proof}
		Note that $$\text{rank} f=\dim \m{Im} f=\dim A-\dim \ker f=\dim A-\text{rank } h,$$hence we have $\dim A=\text{rank} f+\text{rank } h$, similarly, we have $\dim B=\text{rank} g+\text{rank } f$ and $\dim C=\text{rank} h+\text{rank } g$, then it is easy to see we have the desired result.
	\end{proof}

	Now apply the above lemma to our triangle \eqref{trianglecable},  we have 
	
	\begin{equation}\label{FAnk}
		\text{rank}\widehat{F}_{\mathcal{A}_n^{(k)}}=\dfrac{\dim\widehat{\hfk}(S^3,K_{k,kn+1})+\dim \widehat{\hfk}(S^3_n(K),K_n)- \dim\widehat{\hfk}(S^3_{n+1}(K),(K_{n+1})_{k-1,k})}{2},
	\end{equation}where we use the fact that reversing the orientation of a knot preserves the total dimension of its knot Floer homology.

	\subsection{Non-vanishing of $\fl(X_n(K),K_n,\alpha_{-k})$ for $k=0,1,2$}
	
	Let  $-F_{k,kn+1}$ be the  Seifert surface  of the  cable knot $-K_{k,kn+1}$, it is easy to see the surface $A_n^{k}\cup(- F_{k,kn+1})$ represents $\alpha_{-k}=[D^{cc}]-k[\widehat{F}]\in \text{H}_2(X_n(K),K_n)$. Hence for pair $(X_n(K),K_n)$, we can choose a Floer lasagna filling $$\mathcal{F}_{k,v}=\left(\{B,-K_{k,kn+1},v\},(A_n^{(k)},\mathcal{A})\right),$$where $v\in \widehat{\hfk}(S^3,-K_{k,kn+1})$. Then we have $[\mathcal{F}_{k,v}]\in \fl(X_n(K),K_n,\alpha_{-k})$ and $$\m{ev}\left([\mathcal{F}_{k,v}]\right)=\widehat{F}_{\mathcal{A}_n^{(k)}}(v)\in \widehat{\hfk}(S^3_n(K),K_n).$$From the equation above, we can easily deduce that \begin{equation}\label{evineuqlaity}
		\dim \fl(X_n(K),K_n,\alpha_{-k})\ge \m{rank }\ \m{ev}\ge \m{rank }\ \widehat{F}_{\mathcal{A}_n^{(k)}}.
	\end{equation}

	Now to estimate the dimension of $\fl(X_n(K),K_n,\alpha_{-k})$,  it suffices to estimate $ \m{rank }\ \widehat{F}_{\mathcal{A}_n^{(k)}}$. Note that the $(\pm1,q)$-cables of knot $K$ are isotopic to $\pm K$, and $(0,\pm1)$-cable of a knot is  an unknot. Hence we can  apply \eqref{FAnk} to the special case of $k=0,1,2$, 
	\begin{equation}\label{FA0}
		\text{rank}\ \widehat{F}_{\mathcal{A}_n^{(0)}}=\dfrac{\dim\widehat{\hfk}(S^3,U)+\dim \widehat{\hfk}(S^3_n(K),K_n)- \dim\widehat{\hfk}(S^3_{n+1}(K),K_{n+1})}{2},
	\end{equation}\begin{equation}\label{FA1}
		\text{rank}\ \widehat{F}_{\mathcal{A}_n^{(1)}}=\dfrac{\dim\widehat{\hfk}(S^3,K)+\dim \widehat{\hfk}(S^3_n(K),K_n)- \dim\widehat{\hfk}(S^3_{n+1}(K),U)}{2}.
	\end{equation}and\begin{equation}\label{FA2}
		\text{rank}\ 
		\widehat{F}_{\mathcal{A}_n^{(2)}}=\dfrac{\dim\widehat{\hfk}(S^3,K_{2,2n+1})+\dim \widehat{\hfk}(S^3_n(K),K_n)- \dim\widehat{\hfk}(S^3_{n+1}(K),K_{n+1})}{2},
	\end{equation}

	In \cite[Theorems 1.2 and 3.6]{Hed07}, Hedden identified  the knot Floer homology of $K_n\subset S^3_n(K)$ with the knot Floer homology of Whitehead double \(D_+(K,n)\subset S^3\) at Alexander grading \(1\), and computed the latter groups. We can reformulate his result as follows:

	\begin{theorem}[\cite{Hed07}]
		Let $K\subset S^3$ be a knot with Seifert genus $g$, then for $n\ge 2\tau(K)$, we have \begin{equation}\label{ndayutauk}
			\dim \widehat{\hfk}(S^3_n(K),K_n)= n-2g-2+ 2\sum_{a=-g}^g \dim \m{H}\left(C\{i=0,j\le a\}\right).
		\end{equation}And for $n<2\tau (K)$, we have \begin{equation}\label{nxiaoyutauk}
			\dim \widehat{\hfk}(S^3_n(K),K_n)= 4\tau(K)-n-2g-2+ 2\sum_{a=-g}^g \dim \m{H}\left(C\{i=0,j\le a\}\right).
		\end{equation}
	\end{theorem}
	\begin{figure}[htbp]
		\centering
		\includegraphics[width=0.4\linewidth]{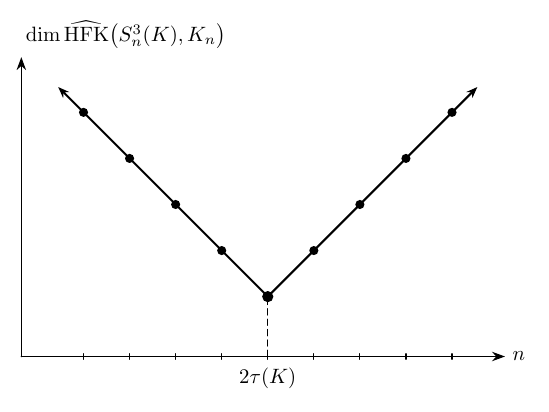}
		\caption{V-shape figure of $\dim \widehat{\hfk}(S^3_n(K),K_n)$}
		\label{fig:dualknotdimensionvshape}
	\end{figure}

	From the theorem above, we obtain  $$\dim \widehat{\hfk}(S^3_n(K),K_n)- \dim\widehat{\hfk}(S^3_{n+1}(K),K_{n+1})=\begin{cases}
		1 & \text{ if } n<2\tau(K) \\
		-1 & \text{ if } n\ge 2\tau(K)
	\end{cases},$$thus from \eqref{FA0} we have 
	\begin{align*}
		\m{rank} \ \widehat{F}_{\mathcal{A}_n^{(0)}}=\dfrac{1+\dim \widehat{\hfk}(S^3_n(K),K_n)- \dim\widehat{\hfk}(S^3_{n+1}(K),K_{n+1})}{2}=\begin{cases}
			1 & \text{ if } n<2\tau(K) \\
			0 & \text{ if } n\ge 2\tau(K)
		\end{cases},
	\end{align*}where we use $\dim \widehat{\hfk}(S^3,U)=1$. Hence from \eqref{evineuqlaity} we can recover the nonvanishing conclusion of Theorem  \ref{thm1.2}, i.e., when $n\le 2\tau(K)-1$, $\fl(X_n(K),K_n,\alpha_0)\ne 0$.

	Now note that when $n\ne -1$, we have $K_{n+1}$ is rationally null-homologous. Hence there is a spectral sequence from $\widehat{\hfk}(S^3_{n+1}(K),K_{n+1})$ to $\widehat{\m{HF}}(S^3_{n+1}(K))$ \cite{OS11}, thus $$\dim\widehat{\hfk}(S^3_{n+1}(K),K_{n+1})\ge \dim\widehat{\m{HF}}(S^3_{n+1}(K))=\dim\widehat{\hfk}(S^3_{n+1}(K),U).$$Combining this fact and \eqref{FA1}, we have for any $n\ne -1$,\begin{align*}
		\text{rank}\ \widehat{F}_{\mathcal{A}_n^{(1)}}&=\dfrac{\dim\widehat{\hfk}(S^3,K)+\dim \widehat{\hfk}(S^3_n(K),K_n)- \dim\widehat{\hfk}(S^3_{n+1}(K),U)}{2}\\
		&\ge \dfrac{\dim\widehat{\hfk}(S^3,K)+\dim \widehat{\hfk}(S^3_n(K),K_n)- \dim\widehat{\hfk}(S^3_{n+1}(K),K_{n+1})}{2}\\
		&\ge \dfrac{\dim\widehat{\hfk}(S^3,K)-1}{2}.
	\end{align*}

	When $n=-1$, from the surgery exact triangle of Heegaard Floer homology, consider the surgery triad $(S^3,S^3_{-1}(K),S^3_0(K))$, and from Lemma \ref{alglem}, we have \begin{align*}
		\dim \widehat{\m{HF}}(S^3_0(K))&\le \dim \widehat{\m{HF}}(S^3)+\dim \widehat{\m{HF}}(S^3_{-1}(K)) \le 1+\dim \widehat{\hfk}(S^3_{-1}(K),K_{-1}).
	\end{align*}

	Hence we again have \begin{align*}
		\text{rank}\ \widehat{F}_{\mathcal{A}_{-1}^{(1)}}&=\dfrac{\dim\widehat{\hfk}(S^3,K)+\dim \widehat{\hfk}(S^3_{-1}(K),K_{-1})- \dim\widehat{\hfk}(S^3_{0}(K),U)}{2}\\
		&\ge \dfrac{\dim\widehat{\hfk}(S^3,K)+\dim \widehat{\hfk}(S^3_{-1}(K),K_{-1})- \left(1+\dim \widehat{\hfk}(S^3_{-1}(K),K_{-1})\right)}{2}\\
		&= \dfrac{\dim\widehat{\hfk}(S^3,K)-1}{2}.
	\end{align*} 
	
	Thus for any integer $n$, we have $$\text{rank}\ \widehat{F}_{\mathcal{A}_n^{(1)}}\ge \dfrac{\dim\widehat{\hfk}(S^3,K)-1}{2}.$$	Now combine this inequality and \eqref{evineuqlaity}, we have $$\dim \mathcal{FL}(X_n(K),K_n,\alpha_{-1})\ge \dfrac{\dim \widehat{\text{HFK}}(K)-1}{2},$$which finishes the proof of Theorem \ref{thm1.1} (1).

	Similarly, one can show that  $$\dim \mathcal{FL}(X_n(K),K_n,\alpha_{-2})\ge \text{rank} \ \widehat{F}_{\mathcal{A}_n^{(2)}} \ge \dfrac{\dim \widehat{\text{HFK}}(K_{2,2n+1})-1}{2},$$which finishes the proof of Theorem \ref{thm1.1} (2).

	\subsection{Calculation for  L-space knots}\label{calculation L space}
	In this section, we  calculate the dimension of $$\widehat{\hfk}(S^3,K_{k,kn+1}), \widehat{\hfk}(S^3_n(K),K_n) \text{ and } \widehat{\hfk}(S^3_{n+1}(K),(K_{n+1})_{k-1,k})$$ when $K$ is an  L-space knot and $n\ge \max\{2g(K),1\}$, $k\ge 2$. Together with \eqref{evineuqlaity}, these computations will prove   Theorem \ref{thm1.5} .

	\begin{definition}
		A rational homology \(3\)-sphere \(Y\) is called an \textit{L-space} if $\dim\widehat{\m{HF}}(Y)=\left| \m{H}_1(Y;\mathbb{Z})\right|$.
	\end{definition}
	
	\begin{definition}
		Suppose that $K\subset Y$ is a framed knot in an L-space, and some positive surgery on $Y$ along $K$ is also an L-space, then we call $K$  an \textit{L-space knot}.
	\end{definition}
	
	Now consider a connected  $3$-manifold $M$ with a single torus boundary component. Let
	\[
	\m{Sl}(M) = \{\theta \in \m{H}_1(\partial M) \mid \theta \text{ is primitive}\}/ \{\pm 1\}
	\]
	be the set of slopes on \(\partial M\). \(\m{Sl}(M)\) is a one-dimensional projective space defined over the rational numbers. If we fix a basis \(\langle \mu, \lambda \rangle\) for \(\m{H}_1(\partial M)\), we can identify \(\m{Sl}(M)\) with \(\overline{\mathbb{Q}} := \mathbb{Q} \cup \{\infty\}\) via the map \(a\mu + b\lambda \mapsto a/b\). We denote by \(M(\theta)\) the closed manifold obtained by Dehn filling \(M\) with slope \(\theta\), and let \(K_\theta \subset M(\theta)\) be the core of the filling solid torus. For example, when $M=S^3\setminus\nu(K)$, and $\theta=p\mu+q\lambda$ as the notation in Section \ref{surgerytraid}, we have $M(\theta)=S^3_{p/q}(K)$ and $K_\theta$ is the dual knot $K_{p/q}$.
	
	\begin{definition}
		If \(M\) is a compact connected oriented $3$-manifold with torus boundary,
		\[
		\mathcal{L}(M) = \{\theta \in \m{Sl}(M) \mid M(\theta) \text{ is an L-space}\}
		\]
		is the set of L-space filling slopes of \(M\). Under the identification of $\m{Sl}(M)$ and $\overline{\mathbb{Q}}$, we denote the interior of $\mathcal{L}(M)$ as $\mathcal{L}^\circ (M)$.
	\end{definition}
	
	\begin{example}\label{2g-1}
		If $M=S^3\setminus \nu(K)$ and $K$ is a nontrivial  L-space knot, then from \cite[Proposition 9.6]{OS11} one can easily deduce  that for any rational number $\dfrac{p}{q}\ge 2g(K)-1$, $S^3_{p/q}(K)$ is an L-space. Hence any rational number $\dfrac{p}{q}\ge 2g(K)-1$ lies in $\mathcal{L}\left(S^3\setminus \nu(K)\right)$. And it is easy to see for any $\dfrac{p}{q}>2g(K)-1$, it lies in the interior $\mathcal{L}^\circ(S^3\setminus\nu(K))$.
	\end{example}

	\begin{definition}
		A knot $K$ in a rational homology sphere $Y$ is called \textit{Floer simple} if
		\[
		\dim\widehat{\hfk}(Y,K)=\dim \widehat{\m{HF}}(Y)=|H_1(Y;\mathbb Z)|.
		\]
	\end{definition}

	Note that J. Rasmussen and S. Rasmussen proved the following criterion for the dual knot to be Floer simple,
	
	\begin{proposition}[{\cite[Corollary 3.6]{RR15}}]
		Let $M$ be a compact connected oriented $3$-manifold with torus boundary, then the dual knot $K_\theta$ in $M(\theta)$ is Floer simple if and only if  $\theta\in \mathcal{L}^\circ (M)$.
	\end{proposition}

	Now combine the above proposition and Example \ref{2g-1}, we have the following corollary,
	\begin{corollary}\label{S3dualknot}
		Let $K\subset S^3$  be a nontrivial  L-space knot, then for any rational number $\dfrac{p}{q}>2g(K)-1$, the dual knot $K_{p/q}$ is Floer simple  in $S^3_{p/q}(K)$. Hence we have $$\dim \widehat{\hfk}\left(S^3_{p/q}(K),K_{p/q}\right)=\left|\m{H}_1\left(S^3_{p/q}(K)\right)\right|=|p|.$$
	\end{corollary}
	
	Then for any L-space knot $K$ and  any integer $n\ge \max\{2g(K),1\}$,  we  have \begin{equation}\label{ndualknot}
		\dim \widehat{\hfk}(S^3_n(K),K_n)=n.
	\end{equation} 	Note that  Hedden also  showed this result in \cite[Theorem 1.4]{Hed11}.

	We now use exact triangles similar to \eqref{trianglecable} to give an upper bound of $\dim \widehat{\hfk}(S^3_{n+1}(K),(K_{n+1})_{k-1,k})$. Recall we take   $J_n^{(k)}=-(kn+1)\mu-k\lambda \subset S^3\setminus\nu(K)$, and $-(K_{n+1})_{k-1,k}$ is the image of $J_n^{(k)}$ in $S^3_{n+1}(K)$.
	
	We now take slope $\theta=n\mu+\lambda$ and a sequence of slopes $\theta_j=(jn+1)\mu+j\lambda$ of $M=S^3\setminus\nu(K)$ for any integer $j\ge 1$. Easy to see $$M(\theta)=S^3_{n}(K),\quad M(\theta_j)=S^3_{n+1/j}(K).$$And as we showed in Section \ref{surgerytraid}, $J_n^{(k)}$ is isotopic to $K_n$ in $M(\theta)$ and isotopic to $-(K_{n+1})_{k-1,k}$ in $M(\theta_1)$. Hence from \eqref{ndualknot} we know that $\dim \widehat{\hfk}\left(M(\theta),J_n^{(k)}\right)=n$, and now we want to determine \begin{equation}\label{5.13}
		\dim \widehat{\hfk}(S^3_{n+1}(K),(K_{n+1})_{k-1,k})=	\dim \widehat{\hfk}\left(M(\theta_1),J_n^{(k)}\right)
	\end{equation}
	
	Note that $\theta_{j+1}=\theta_j+\theta$ and $\theta\cdot \theta_j=1$, then $$\theta_j\cdot \theta=\theta\cdot (-\theta_{j+1})=(-\theta_{j+1})\cdot \theta_j=-1,$$ hence for the triple $(\theta_j,\theta,-\theta_{j+1})$,  we have an exact triangle followed from Theorem \ref{knotexact} for any integer $j\ge 1$, \begin{equation}\label{jtriangle}
		\begin{tikzcd}
			{\widehat{\mathrm{HFK}}\left(M(\theta_j),J_n^{(k)}\right)} && {\widehat{\mathrm{HFK}}\left(M(\theta),J_n^{(k)}\right)} \\
			& {\widehat{\mathrm{HFK}}\left(M(-\theta_{j+1}),J_n^{(k)}\right)}
			\arrow[from=1-1, to=1-3]
			\arrow[from=1-3, to=2-2]
			\arrow[from=2-2, to=1-1]
		\end{tikzcd}
	\end{equation}

	Since $M(-\theta_{j+1})$ and  $M(\theta_{j+1})$ denote the same manifold, then we have $$\dim \widehat{\mathrm{HFK}}\left(M(-\theta_{j+1}),J_n^{(k)}\right)=\dim \widehat{\mathrm{HFK}}\left(M(\theta_{j+1}),J_n^{(k)}\right).$$ It follows from this exact triangle and Lemma \ref{alglem} that \begin{equation*}
		\dim \widehat{\mathrm{HFK}}\left(M(\theta_j),J_n^{(k)}\right)-\dim \widehat{\mathrm{HFK}}\left(M(\theta_{j+1}),J_n^{(k)}\right)\le \dim \widehat{\mathrm{HFK}}\left(M(\theta),J_n^{(k)}\right)=n.
	\end{equation*}Iterating this inequality yields \begin{equation}\label{5.15}
		\dim \widehat{\mathrm{HFK}}\left(M(\theta_1),J_n^{(k)}\right)\le (k-2)n+ \dim \widehat{\mathrm{HFK}}\left(M(\theta_{k-1}),J_n^{(k)}\right).
	\end{equation}

	Note that in $M(\theta_{k-1})$, the Dehn filling slope is $\theta_{k-1}=\left((k-1)n+1\right)\mu+(k-1)\lambda$, and $J_n^{(k)}=-(kn+1)\mu-k\lambda=-\theta_k$, hence $J_n^{(k)}\cdot \theta_{k-1}=(-\theta_k)\cdot \theta_{k-1}=-1$. Then we know that $J_n^{(k)}$ is isotopic to the core of this filling torus, i.e., the dual knot with orientation reversed $-K_{\theta_{k-1}}$, hence \begin{equation}\label{5.16}
		\dim \widehat{\mathrm{HFK}}\left(M(\theta_{k-1}),J_n^{(k)}\right)=\dim   \widehat{\mathrm{HFK}}\left(M(\theta_{k-1}),K_{\theta_{k-1}}\right).
	\end{equation}
	
	Since $\theta_{k-1}$ corresponds to $\dfrac{(k-1)n+1}{k-1}=n+\dfrac{1}{k-1}$, now from $n\ge \max\{2g(K),1\}$ and $k\ge 2$, we have the $n+\dfrac{1}{k-1}>2g(K)-1$. Hence from the Corollary  \ref{S3dualknot}, we have $K_{\theta_{k-1}}$ is Floer simple, and 
	\begin{equation}\label{5.17}
		\dim\widehat{\mathrm{HFK}}\left(M(\theta_{k-1}),K_{\theta_{k-1}}\right)=(k-1)n+1.
	\end{equation}

	Now combine \eqref{5.13}, \eqref{5.15},\eqref{5.16} and \eqref{5.17}, we have \begin{equation}
		\dim \widehat{\hfk}(S^3_{n+1}(K),(K_{n+1})_{k-1,k})\le (2k-3)n+1.
	\end{equation}

	It remains to determine the dimension of $\widehat{\hfk}(S^3,K_{k,kn+1})$, note that 
	
	\begin{theorem}[\cite{Hom15}]
		If $K\subset S^3$ is an oriented knot, then $K_{p,q}$ is an L-space knot if and only if $K$ is an  L-space knot and $\dfrac{q}{p}\ge 2g(K)-1$, where $p>1$, $q>0$ and $\gcd(p,q)=1$.
	\end{theorem}
	
	Now for our case, $n\ge 2g(K)$ and $k\ge 2$, hence $\dfrac{kn+1}{k}=n+\dfrac{1}{k}\ge 2g(K)-1$, then the cabling criterion therefore implies that \(K_{k,kn+1}\) is an L-space knot. 
	
	From \cite[Theorem 1.2 and Corollary 1.3]{OS05} ,  for the L-space knots in $S^3$,  their knot Floer complexes are staircases and can be determined from their Alexander polynomials. More precisely, let $\Delta_K(t)$ be the Alexander polynomial of $K$, 
	then its nonzero coefficients  are all $\pm 1$ and alternate in sign, beginning and ending with $+1$, and $$\dim \widehat{\hfk}(K)=\left\| \Delta_K(t)\right\|_1,$$where $\left\|\cdot\right\|_1$ is a norm on $\mathbb{Z}[t,t^{-1}]$ determined by $\displaystyle \left\|\sum_{i=-n}^m a_it^i\right\|_1:=\sum_{i=-n}^{m}|a_i|$.

	Now from  $K_{k,kn+1}$ is an L-space knot, we have \begin{equation}\label{cablealexander}
		\dim \widehat{\hfk}\left( K_{k,kn+1}\right)=\left\| \Delta_{K_{k,kn+1}}(t)\right\|_1,
	\end{equation}hence it suffices to determine the Alexander polynomial of $K_{k,kn+1}$. 	From the cabling formula of Alexander polynomial, we have 
	\begin{align*}
		\Delta_{K_{k,kn+1}}(t)&=\Delta_K(t^k)\cdot \Delta_{T_{k,kn+1}}(t)\\
		&=\Delta_K(t^k)\cdot \dfrac{(t^{k(kn+1)}-1)(t-1)}{(t^k-1)(t^{kn+1}-1)}\\
		&=(1-t)\cdot \dfrac{\Delta_K(t^k)}{1-t^k}\cdot \left(\sum_{j=0}^{k-1}(t^{kn+1})^j\right)\\
		&=(1-t)\cdot H(t^k)\cdot \left(\sum_{j=0}^{k-1}(t^{kn+1})^j\right),
	\end{align*}where we denote $H(t):=\dfrac{\Delta_K(t)}{1-t}$ for short. Every exponent occurring in \(H(t^k)\) is divisible by \(k\). Collecting terms according to their residue classes modulo \(k\) gives a decomposition \begin{equation}
		\Delta_{K_{k,kn+1}}(t)=\sum_{j=0}^{k-1} t^j P_j(t^k).
	\end{equation}And it is easy to see \begin{equation}
		\left\|\Delta_{K_{k,kn+1}}(t)\right\|_1=\sum_{j=0}^{k-1}\left\| P_j(t^k)\right\|_1.
	\end{equation}

	Now note that for $1\le j\le k-1$, $$P_j(t^k)=(t^{kjn}-t^{k(j-1)n})\cdot H(t^k)= t^{k(j-1)n}(t^{kn}-1)H(t^k),$$hence we have $\left\| P_j(t^k)\right\|_1=\left\| (t^{kn}-1)H(t^k)\right\|_1=\left\|  (1-t^n)H(t)\right\|_1$.
	
	For $j=0$, we have $$P_0(t^k)=(1-t^{k((k-1)n+1)})\cdot H(t^k),$$hence $\left\| P_0(t^k)\right\|_1=\left\|  (1-t^{(k-1)n+1})H(t)\right\|_1$. Hence if we denote $\rho(N)$ as $\left\|  (1-t^N)H(t)\right\|_1$, then  from \eqref{cablealexander} we have
	\begin{equation}\label{kn+1}
		\dim \widehat{\hfk}(K_{k,kn+1})=(k-1)\cdot \rho(n)+\rho\left((k-1)n+1\right).
	\end{equation}
	
	Now it remains to compute $\rho(N)=\left\|  (1-t^N)H(t)\right\|_1$. We claim that, whenever $N\geq\max\{2g(K),1\}$,
	\begin{equation}\label{rhoN}
		\rho(N)=\left\|  (1-t^N)H(t)\right\|_1=N.
	\end{equation}

	Recall $$(1-t^N)H(t)=(1-t^N)\cdot \dfrac{\Delta_K(t)}{1-t}=(1+t+\cdots+t^{N-1})\Delta_K(t).$$ Since $K$ is L-space knot, we can assume $\Delta_K(t)=\displaystyle\sum_{i=0}^{2m} (-1)^i t^{a_i}$, such that $a_0=0$ and $a_{2m}=2g(K)$.
	
	Let $c_j$ be the sum of the coefficients of $\Delta_K(t)$ in degrees
	at most $j$, with $c_j=0$ for $j<0$.
	Alternation gives $c_j\in\{0,1\}$, and $c_j=1$ for
	$j\geq\deg \Delta_K(t)=2g(K)$. Easy to see 
	the coefficient in degree $j$ of $(1-t^N)H(t)=(1+t+\cdots+t^{N-1})\Delta_K(t)$ is $c_j-c_{j-N}$.
	For $j<N$ the second term is zero, while for $j\geq N$
	the first term is one. Thus every coefficient of  $(1-t^N)H(t)$ is zero or one. Evaluating at \(t=1\) and using \(\Delta_K(1)=1\) gives $\left\|  (1-t^N)H(t)\right\|_1=N$. Hence we know $\rho(N)=N$.

	Now combine \eqref{kn+1} and \eqref{rhoN}, we have \begin{equation}\label{dimcable}
		\dim \widehat{\hfk}(K_{k,kn+1})=2(k-1)n+1.
	\end{equation}

	In summary, we have \begin{align*}
		\dim \widehat{\hfk}(K_{k,kn+1})&=2(k-1)n+1,\\
		\dim \widehat{\hfk}(S^3_n(K),K_n)&=n,\\
		\dim \widehat{\hfk}(S^3_{n+1}(K),(K_{n+1})_{k-1,k})&\le (2k-3)n+1.
	\end{align*}
	
	From the exact triangle \eqref{trianglecable} and Lemma \ref{alglem}, we have \begin{align*}
		0&\le 2\cdot \m{rank }\ \widehat{F}_2\\
		&= \dim \widehat{\hfk}(S^3_n(K),K_n)+\dim \widehat{\hfk}(S^3_{n+1}(K),(K_{n+1})_{k-1,k})	\\
		&\quad -\dim \widehat{\hfk}(K_{k,kn+1})\\
		&\le n+(2k-3)n+1-\left(2(k-1)n+1\right)=0.
	\end{align*} Hence equality holds throughout. In particular, we have \begin{equation}
		\dim \widehat{\hfk}(S^3_{n+1}(K),(K_{n+1})_{k-1,k})= (2k-3)n+1.
	\end{equation}

	Finally, from \eqref{FAnk}, we have \begin{align*}
		\text{rank}\widehat{F}_{\mathcal{A}_n^{(k)}}&=\dfrac{\dim\widehat{\hfk}(S^3,K_{k,kn+1})+\dim \widehat{\hfk}(S^3_n(K),K_n)- \dim\widehat{\hfk}(S^3_{n+1}(K),(K_{n+1})_{k-1,k})}{2}\\
		&=\dfrac{2(k-1)n+1+n-\left((2k-3)n+1\right)}{2}=n.
	\end{align*}Hence using \eqref{evineuqlaity}, we have $$\dim \fl \left(X_n(K),K_n,\alpha_{-k}\right)\ge \text{rank }\widehat{F}_{\mathcal{A}_n^{(k)}}=n,$$which finishes the proof of Theorem \ref{thm1.5}.

	\section{Further Questions}\label{sec6}

	In this section we propose some further questions about the Floer lasagna modules.
	
	First of all, it is natural to compare the theory of Floer lasagna modules and skein lasagna modules. In \cite[Theorem 1.4]{MN22}, Manolescu and Neithalath proved a connected sum formula for skein lasagna modules, hence we can ask whether we can prove a similar result for Floer lasagna modules.
	
	\begin{question}
		Is there an  isomorphism between $$\fl(W_1\natural W_2,\mathbb{L}_1\# \mathbb{L}_2)\quad \text{and}\quad \fl(W_1,\mathbb{L}_1)\otimes \fl(W_2,\mathbb{L}_2)\ ? $$
	\end{question}
	
	\begin{remark}
		Note that in \cite[Theorem 1.4]{MN22} they use the disjoint union of links rather than connected sum. We cannot use disjoint union since we have $\dim \widehat{\hfl}(S^3,U\cup U)=2$ is greater than $\dim \left(\widehat{\hfl}(S^3,U)\otimes \widehat{\hfl}(S^3,U)\right)=1$. Thus the unmodified disjoint-union formula cannot hold in the Floer setting. This motivates considering a connected-sum formula for based links and hence one can not directly copy the proof of  \cite[Theorem 1.4]{MN22}.
	\end{remark}
	
	\begin{remark}
		If we have the Floer lasagna module connected sum isomorphism, then connected sum $(W,U)$ with $(B^4,L)$ where $U$ is local unknot we get $(W,L)$. Then from the isomorphism  $\fl(B^4,L)\cong \widehat{\hfl}(S^3,L)\ne 0$ \cite[Lemma 2.8]{chen2022floer}, we  have $\fl(W,L)\ne 0$ for any local link $L\subset B^3\subset \partial W$ if $\fl(W,U)\ne 0$.
	\end{remark}

	\begin{question} 
		Can the local unknot in Proposition \ref{thmstein} and \ref{thm1.0} be replaced by the empty link? More precisely,
		\begin{itemize}
			\item for every Stein domain \(X\), is \(\fl\left(X\#^n\overline{\mathbb{C}P^2}\right)\) nonzero?
			
			\item for every knot \(K\) and every integer \(n\le \widehat{\nu}(K)+0^{\widehat{\nu}(K)}-1\), is \(\fl(X_n(K))\) nonzero?
		\end{itemize}
	\end{question}
	
	\begin{remark}
		For the question related to Proposition \ref{thm1.0}, easy to see, when $K$ is slice, then $\nu(K)=\varepsilon(K)=0$, hence $\widehat{\nu}(K)=0$. Thus we have shown the  above question has an affirmative answer when $K$ is slice in Proposition \ref{thmknottrace}. For the left-handed trefoil \(3_1\), we have \(\nu(3_1)=0\) and \(\varepsilon(3_1)=-1\) thus $\widehat{\nu}(3_1)=-1$, so the question concerns all \(n\le-2\). Remark \ref{31knot} gives an affirmative answer for \(n\le-4\).
	\end{remark}
	
	\begin{remark}
		Cai Zeng pointed out to the author that, for a general $4$-manifold $W$ with nonempty boundary and $U\subset \partial W$ a local unknot, one can glue $\left(\partial W\times I,(D,\mathcal{A}_D)\right)$ to the Floer lasagna filling of $(W,U)$ to obtain a Floer lasagna filling of $(W,\varnothing)$, and the gluing map induced by this cobordism is a surjection from $\fl(W,U)$ to $\fl(W)$. However, it is not clear to the author whether this gluing map is injective, or even whether it sends $\mathcal{F}_{W,U}$ to a nonzero element. If either holds, then we would have $[\mathcal{F}_W]\ne 0$ and $\fl(W)\ne 0$.
	\end{remark}

The calculation for L-space knots in Section \ref{calculation L space} relies
mainly on the fact that the Alexander polynomial determines the knot Floer
homology of L-space knots. In \cite{Bin23}, Binns classified the knot Floer
complex of almost L-space knots. Consequently, the Alexander polynomial also
determines their knot Floer homology. More recently, Zou \cite{Zou26}
classified the knot Floer complex of height $4$ knots. This extends the cases
of height $0$ and height $2$, which correspond to L-space knots and almost
L-space knots, respectively. Hence it is natural to ask whether the calculation
in Section \ref{calculation L space} can be generalized to knots of height $m$,
where $m$ is an even integer. The author is grateful to Cai Zeng for proposing
the following question.
	
	\begin{question}
		If $K\subset S^3$ is a height $m$  knot, does there exist an integer $N_K$, depending on  $K$, such that $\fl(X_n(K),K_n)$ is infinite-dimensional for every integer $n\ge N_K$?
	\end{question}
	
	The main technical results we proved in this paper are all about the nonvanishing of some knot Floer cobordism maps. For example we constructed an annulus cobordism between $-K\subset S^3$ and $K_n\subset S^3_n(K)$ in Section \ref{surgerytraid}, and proved that this cobordism map induces nonzero knot Floer cobordism map. Hence it is natural to ask,
	
	\begin{question}
		Is it possible to use the annulus cobordism map between $-K$ and $K_n$ to construct exotic surfaces bounded by $K_n$ in $X_n(K)$?
		
		More precisely, one may have exotic surfaces bounded by $-K$ in $4$-ball, and the exoticness can be detected by knot Floer homology. Composing the exotic surfaces with our annulus,  we have two surfaces bounded by $K_n$ in $X_n(K)$ and are topologically isotopic rel boundary. One may use the nonvanishing of the annulus cobordism map to show that after the composition the surfaces  remain  exotic.
	\end{question}

	\bibliographystyle{alpha}

\begin{thebibliography}{EVVZ17}

\bibitem[AGL26]{AGL26}
  Akram Alishahi, Eugene Gorsky, and Beibei Liu.
  \newblock Floer lasagna modules for some 2-handlebodies.
  \newblock {{arXiv preprint: 2610.34744}}, 2026.

\bibitem[Bin23]{Bin23}
  Fraser Binns.
  \newblock The {CFK}\(^{\infty}\) type of almost \(L\)-space knots.
  \newblock {{arXiv preprint: 2303.07249}}, 2023.

\bibitem[Che22]{chen2022floer}
  Daren Chen.
  \newblock Floer lasagna modules from link Floer homology.
 \newblock {\textit{Algebr. Geom. Topol.}}, to appear, 2026.
  \newblock arXiv preprint: 2203.07650.

\bibitem[Eft18]{Eft18}
  Eaman Eftekhary.
  \newblock Bordered {F}loer homology and existence of incompressible tori in homology spheres.
  \newblock {\textit{Compos. Math.}}, 154(6): 1222--1268, 2018.

\bibitem[Eli90]{Eli90}
  Yakov Eliashberg.
  \newblock Filling by holomorphic discs and its applications.
  \newblock In {\textit{Geometry of Low-Dimensional Manifolds, 2 (Durham, 1989)}}, volume 151 of \textit{London Math. Soc. Lecture Note Ser.}, pages 45--67. Cambridge Univ. Press, 1990.

\bibitem[GS99]{GS99}
  Robert  Gompf and András Stipsicz.
  \newblock \textit{4-Manifolds and Kirby Calculus}.
  \newblock Graduate Studies in Mathematics, 20. American Mathematical Society, Providence, RI, 1999.

\bibitem[Hed07]{Hed07}
  Matthew Hedden.
  \newblock Knot Floer homology of Whitehead doubles.
  \newblock {\textit{Geom. Topol.}}, 11(4): 2277--2338, 2007.

\bibitem[Hed08]{Hed08}
  Matthew Hedden.
  \newblock An Ozsv{\'a}th-Szab{\'o} Floer homology invariant of knots in a contact manifold.
  \newblock {\textit{Adv. Math.}}, 219(1): 89--117, 2008.

\bibitem[Hed09]{Hed09}
  Matthew Hedden.
  \newblock On knot Floer homology and cabling. II.
  \newblock {\textit{Int. Math. Res. Not. IMRN}}, 2009(12): 2248--2274, 2009.

\bibitem[Hed11]{Hed11}
  Matthew Hedden.
  \newblock On {F}loer homology and the {B}erge conjecture on knots admitting lens space surgeries.
  \newblock {\textit{Trans. Amer. Math. Soc.}}, 363(2): 949--968, 2011.

\bibitem[HHSZ26]{HHSZ26}
  Kristen Hendricks, Jennifer Hom, Matthew Stoffregen, and Ian Zemke.
  \newblock An involutive dual knot surgery formula.
  \newblock {\textit{J. Symplectic Geom.}}, 24(2): 227--302, 2026.

\bibitem[HL24]{HL24}
  Matthew Hedden and Adam Simon Levine.
  \newblock A surgery formula for knot Floer homology.
  \newblock {\textit{Quantum Topol.}}, 15(2): 229--336, 2024.

\bibitem[HN22]{HN22}
  Matthew Hedden and Lev Tovstopyat-Nelip.
  \newblock On naturality of the Ozsváth–Szabó contact invariant.
  \newblock {\textit{The Open Book Series}}, 5(1): 123--143, 2022.

\bibitem[Hom11]{Hom15}
  Jennifer Hom.
  \newblock A note on cabling and \(L\)-space surgeries.
  \newblock {\textit{Algebr. Geom. Topol.}}, 11(1): 219--223, 2011.

\bibitem[Hom14]{Hom14}
  Jennifer Hom.
  \newblock The knot {F}loer complex and the smooth concordance group.
  \newblock {\textit{Comment. Math. Helv.}}, 89(3): 537--570, 2014.

\bibitem[HP25]{HP19}
  Kyle Hayden and Lisa Piccirillo.
  \newblock The trace embedding lemma and spinelessness.
  \newblock {\textit{J. Differential Geom.}}, 131(1): 133--163, 2025.


\bibitem[HW26]{HW26}
  Jennifer Hom and Shunyu Wan.
  \newblock Heegaard Floer knot trace invariants, exotic 4-manifolds, and symplectic obstructions.
  \newblock arXiv preprint: 2608.09734, 2026.


\bibitem[Juh06]{Juh06}
  András Juhász.
  \newblock Holomorphic discs and sutured manifolds.
  \newblock {\textit{Algebr. Geom. Topol.}}, 6: 1429--1457, 2006.

\bibitem[Juh16]{Juh16}
  András Juhász.
  \newblock Cobordisms of sutured manifolds and the functoriality of link {F}loer homology.
  \newblock {\textit{Adv. Math.}}, 299: 940--1038, 2016.

\bibitem[JZ20]{JZ21}
  András Juhász and Ian Zemke.
  \newblock Contact handles, duality, and sutured {F}loer homology.
  \newblock {\textit{Geom. Topol.}}, 24(1): 179--307, 2020.

\bibitem[JZ24]{JZ24}
  András Juhász and Ian Zemke.
  \newblock Stabilization distance bounds from link Floer homology.
  \newblock {\textit{J. Topol.}}, 17(2): e12338, 2024.



\bibitem[LY24]{LY24}
Zhenkun Li and Fan Ye.
\newblock  An enhanced Euler characteristic of sutured instanton homology.
\newblock  {\textit{Int. Math. Res. Not.}} 2024(4): 2873–2936, 2024



\bibitem[LY25]{LY25}
Zhenkun Li and Fan Ye.
\newblock  Knot surgery formulae for instanton Floer homology, I: The main theorem.
\newblock {\textit{Geom. Topol.}}, 29(5): 2269--2342, 2025.



\bibitem[McC26]{McC26}
  Colin McCulloch.
  \newblock Computations of Floer lasagna modules for traces of negative $L$-space knots.
  \newblock arXiv preprint: 2609.17503.

\bibitem[MN22]{MN22}
  Ciprian Manolescu and Ikshu Neithalath.
  \newblock Skein lasagna modules for 2-handlebodies.
  \newblock {\textit{J. Reine Angew. Math.}}, 788: 37--76, 2022.

\bibitem[MWW22]{MWW22}
  Scott Morrison, Kevin Walker, and Paul Wedrich.
  \newblock Invariants of 4-manifolds from Khovanov-Rozansky link homology.
  \newblock {\textit{Geom. Topol.}}, 26(8): 3367--3420, 2022.

\bibitem[MWW23]{MWW23}
  Ciprian Manolescu, Kevin Walker, and Paul Wedrich.
  \newblock Skein lasagna modules and handle decompositions.
  \newblock {\textit{Adv. Math.}}, 425: Article 109071, 2023.

\bibitem[OS03]{OS03}
  Peter Ozsv{\'a}th and Zolt{\'a}n Szab{\'o}.
  \newblock Absolutely graded Floer homologies and intersection forms for four-manifolds with boundary.
  \newblock {\textit{Adv. Math.}}, 173(2): 179--261, 2003.



\bibitem[OS04a]{OS04a}
  Peter Ozsv{\'a}th and Zolt{\'a}n Szab{\'o}.
  \newblock Holomorphic disks and topological invariants for closed three-manifolds.
  \newblock {\textit{Ann. of Math. (2)}}, 159(3): 1027--1158, 2004.

\bibitem[OS04b]{OS04b}
  Peter Ozsv{\'a}th and Zolt{\'a}n Szab{\'o}.
  \newblock Holomorphic disks and three-manifold invariants: properties and applications.
  \newblock {\textit{Ann. of Math. (2)}}, 159(3): 1159--1245, 2004.

\bibitem[OS04c]{OS04c}
  Peter Ozsv{\'a}th and Zolt{\'a}n Szab{\'o}.
  \newblock Holomorphic disks and knot invariants.
  \newblock {\textit{Adv. Math.}}, 186(1): 58--116, 2004.

\bibitem[OS05]{OS05}
  Peter  Ozsv{\'a}th and Zolt{\'a}n Szab{\'o}.
  \newblock On knot {F}loer homology and lens space surgeries.
  \newblock {\textit{Topology}}, 44(6): 1281--1300, 2005.

\bibitem[OS06]{OS06}
  Peter Ozsv{\'a}th and Zolt{\'a}n Szab{\'o}.
  \newblock Holomorphic triangles and invariants for smooth four-manifolds.
  \newblock {\textit{Adv. Math.}}, 202(2): 326--400, 2006.

\bibitem[OS08a]{OS08}
  Peter  Ozsváth and Zoltán Szabó.
  \newblock Knot Floer homology and integer surgeries.
  \newblock {\textit{Algebr. Geom. Topol.}}, 8(1): 101--153, 2008.

\bibitem[OS08b]{OS08b}
  Peter Ozsv{\'a}th and Zolt{\'a}n Szab{\'o}.
  \newblock Holomorphic disks, link invariants and the multi-variable {A}lexander polynomial.
  \newblock {\textit{Algebr. Geom. Topol.}}, 8(2): 615--692, 2008.

\bibitem[OS11]{OS11}
  Peter  Ozsv{\'a}th and Zolt{\'a}n Szab{\'o}.
  \newblock Knot {F}loer homology and rational surgeries.
  \newblock {\textit{Algebr. Geom. Topol.}}, 11(1): 1--68, 2011.

\bibitem[Pla04]{Pla04}
  Olga Plamenevskaya.
  \newblock Bounds for the Thurston-Bennequin number from Floer homology.
  \newblock {\textit{Algebr. Geom. Topol.}}, 4(1): 399--406, 2004.

\bibitem[RR17]{RR15}
  Jacob Rasmussen and Sarah Dean Rasmussen.
  \newblock Floer simple manifolds and L-space intervals.
  \newblock {\textit{Adv. Math.}}, 322: 738--805, 2017.

\bibitem[RW26]{RW26}
  Qiuyu Ren and Michael Willis.
  \newblock Khovanov homology and exotic $4$-manifolds.
  \newblock {\textit{Ann. of Math.}}, to appear.
  \newblock arXiv:2402.10452.

\bibitem[SZ24]{SZ24}
  Ian  Sullivan and Melissa Zhang.
  \newblock Kirby belts, categorified projectors, and the skein lasagna module of \(S^{2}\times S^{2}\).
  \newblock {\textit{Quantum Topol.}},  2024.

\bibitem[Zem19a]{Zem19a}
  Ian Zemke.
  \newblock Link cobordisms and functoriality in link Floer homology.
  \newblock {\textit{J. Topol.}}, 12(1): 94--220, 2019.

\bibitem[Zem19b]{Zem19b}
  Ian Zemke.
  \newblock Link cobordisms and absolute gradings on link {F}loer homology.
  \newblock {\textit{Quantum Topol.}}, 10(2): 207--323, 2019.

\bibitem[Zem21]{Zem21}
  Ian Zemke.
  \newblock Duality and mapping tori in {H}eegaard {F}loer homology.
  \newblock {\textit{J. Topol.}}, 14(3): 1027--1112, 2021.

\bibitem[Zou26]{Zou26}
  Shengyu Zou.
  \newblock Knot {F}loer complexes of generalized {L}-space knots.
  \newblock {{arXiv preprint: 2610.35305}}, 2026.

\end{thebibliography}

\end{document}